\documentclass{article}
\usepackage{amsmath}
\usepackage{amsthm}
\usepackage{indentfirst}
\usepackage{amsfonts}
\usepackage{amssymb}
\usepackage{graphicx}
\numberwithin{equation}{section}
\newtheorem{thm}{Theorem}[section]
\newtheorem{lem}{Lemma}[section]
\newtheorem{dfn}{Definition}[section]
\newtheorem{rem}{Remark}[section]
\newtheorem{prop}{Proposition}[section]
\newtheorem{exam}{Example}[section]

\newcommand{\fff}{{\mathbf{f}}}

\begin{document}

\title{Coupling Tensors and Poisson Geometry Near a Single Symplectic Leaf}
\author{Yurii Vorobjev\thanks{Research partially supported by CONACYT Grants 28291-E
and 35212-E.}\\{\small Department of Applied Mathematics,}\\{\small Moscow State Institute of Electronics and Mathematics}\\{\small Moscow, Russia,109028.}\\{\small E-mail address: vorob@amath.msk.ru}\\\medskip\\{\small Departamento de Matem\'{a}ticas}\\{\small Universidad de Sonora}\\{\small Hermosillo, M\'{e}xico, 83000.}\\{\small Current e-mail address:yurimv@guaymas.uson.mx}}
\date{}
\maketitle
\begin{abstract}
In the framework of the connection theory, a contravariant analog of the
Sternberg coupling procedure is developed for studying a natural class of
Poisson structures on fiber bundles, called coupling tensors. We show that
every Poisson structure near a closed symplectic leaf can be realized as a
coupling tensor. Our main result is a geometric criterion for the neighborhood
equivalence between Poisson structures over the same leaf. This criterion
gives a Poisson analog of the relative Darboux theorem due to Weinstein.
Within the category of the algebroids, coupling tensors are introduced on the
dual of the isotropy of a transitive Lie algebroid over a symplectic base. As
a basic application of these results, we show that there is a well defined
notion of a ``linearized'' Poisson structure over a symplectic leaf which
gives rise to a natural model for the linearization problem. \medskip

\textit{MC classification}: 58F05, 53C30, 58H05.

\end{abstract}

\section{Introduction}

The notion of a coupling form due to Sternberg \cite{St} naturally arises from
the study of fiber compatible (pre)symplectic structures on the total space of
a symplectic bundle (for various aspects of this problem, see, for example,
\cite{Tu,We2,We3,GoLSW,GSt,GLS,KV1}). This notion is based on the concept of
connection and curvature and can be introduced for a wide class of bundles
\cite{GLS}. A derivation of the coupling procedure for the associated bundle
$P\times_{G}\frak{g}^{\ast}$ (called the universal phase space) via reduction
was suggested in \cite{We2}.

We are interested in a contravariant version of the Sternberg coupling
procedure in the Poisson category. Suppose we start with a (locally trivial)
Poisson fiber bundle $(E\xrightarrow{\pi}B,{\mathcal{V}}^{\text{fib}})$
equipped with a smooth field ${\mathcal{V}}^{\text{fib}}=\{B\ni{b}%
\mapsto{\mathcal{V}}_{b}^{\text{fib}}\in\chi^{2}(E_{b})\}$ of Poisson
structures on the fibers. Unlike the symplectic case \cite{GoLSW}, the
fiberwise Poisson structure ${\mathcal{V}}^{\text{fib}}$ admits always a
unique extension to a vertical Poisson tensor ${\mathcal{V}}$ on the total
space $E$. Every Ehresmann connection $\operatorname{Hor}\rightarrow
TE\xrightarrow{\Gamma}\operatorname{Vert}$ gives rise to the space of
horizontal multivector fields on $E$. A connection is called Poisson if the
parallel transport preserves ${\mathcal{V}}^{\text{fib}}$. Given a Poisson
connection $\Gamma$, the problem is to describe Poisson bivector fields on $E$
of the form: $\Pi=(\text{$\Gamma$-horizontal bivector field})+{\mathcal{V}}$.
Putting this decomposition into the Jacobi identity, we get two quadratic
equations for the horizontal part of $\Pi$: (I) the horizontal Jacobi identity
and (II) the curvature identity. Under the assumption: $\Pi$ is nondegenerate
on the annihilator of the vertical subbundle $\operatorname{Vert}$, equations
(I) and (II) are reduced to linear equations for a horizontal $2$-form
$\mathbb{F}$. If $\mathbb{F}$ is a solution of these equations, then the
corresponding Poisson tensor $\Pi$ is just what we call a coupling tensor
associated with data $(\Gamma,{\mathcal{V}},\mathbb{F})$. Note that in the
case of coadjoint bundles, Poisson structures of such a type arise naturally
from the study of Wong's equations \cite{MoMR,Mo}.

In this paper we give a systematic treatment of coupling tensors. The first
our observation is that: in a tubular neighborhood $E$ of a closed symplectic
leaf $B$ every Poisson structure $\Psi$ is realized as a coupling tensor
(locally, this follows from the splitting theorem \cite{We4}). As a
consequence, $\Psi$ induces an intrinsic Poisson connection $\Gamma$ on $E$
which gives rise to a geometric characteristic of the leaf. Moreover, the
vertical part ${\mathcal{V}}$ of the coupling tensor at the leaf $B$ is of
rank $0$ and gives a ``global'' realization of local transverse Poisson
structures \cite{We4}. If the symplectic leaf $B$ is regular, then
$\mathcal{V}=0$ and the coupling tensor is the horizontal lift of the
nondegenerate Poisson structure on $B$ via the flat connection associated with
the symplectic foliation. We are interested in the non-flat case when the rank
of a Poisson structure $\Psi$ is not locally constant at $B$.

So, for the study of Poisson structures near a single symplectic leaf we may
restrict our attention to the class of coupling tensors. Our main result is a
geometric criterion for a neighborhood equivalence of two coupling near a
common closed symplectic leaf. The criterion is formulated in terms of the
intrinsic Poisson connection and its curvature and implies a Poisson version
of the relative Darboux theorem \cite{We1}. This result continues our previous
investigations of the formal Poisson equivalence \cite{IKV} . To state the
result, we use a contravariant analog of the homotopy method due to Moser
\cite{Mos} and Weinstein \cite{We1}. The technical part is based on the
Schouten calculus \cite{Li,KM,KV2,Va,K-SM} and the (covariant) connection
theory for general fiber bundles \cite{GHV,GLS}. Note also that a geometric
approach, based on the notion of a contravariant Poisson connection due to
Vaisman \cite{Va}, was developed in \cite{Fe}.

As a basic application of the Poisson neighborhood theorem, we show that there
exists a well defined notion of a linearized Poisson structure over a closed
symplectic leaf which is well known in the zero-dimensional case \cite{We4}.
The linearized Poisson structure is completely determined by the transitive
Lie algebroid of a symplectic leaf \cite{KV2}. To derive this fact, we
introduce and study a class of coupling tensors associated with transitive Lie
algebroids over a symplectic base. This class consists of isomorphic Poisson
structures parametrized by connections on the Lie algebroid in the sense of
Mackenzie \cite{Mz}. Here we use\ the technique of adjoint connections on the
dual of a Lie algebriod. Adjoint connections naturally arise in the general
theory of Lie algebriods \cite{Mz} (also see \cite{Ku} ) as well as in the
context of infinitesimal Poisson geometry \cite{KV2,IKV}. We show that the
holonomy of adjoint connections is related with the notion of a linear Poisson
holonomy introduced in \cite{GiGo} (the definition of this notion in terms of
contravariant connections can be found in \cite{Fe} ).

In this paper we do not discuss the linearization problem. But we hope that
the linearized Poisson model, introduced here, can be used for extension of
results on the linearizability at a point \cite{We4,Cn1,Cn2,Du} to the
higher-dimentional case.

The body of the paper is organized as follows. In Section~2, a general
description of coupling tensors in terms of geometric data is given in
Theorem~2.1. In Section~3, we formulate main results on a neighborhood
equivalence of two coupling tensors with the same symplectic leaf,
Theorem~3.1. and Theorem 3.2. The important technical part of the proof of
Lemma~3.1., is given in Appendix~A. Section~4 is devoted to coupling tensors
associated with transitive Lie algebroids. Here we show that the criterion in
Theorem~3.1 leads to the equivalence relation for Lie algebroids. In
Section~5, using results of Section~4, we give a definition of the linearized
Poisson structure of a symplectic leaf and discuss some applications.

\subsection*{Acknowledgments}

The main results of the paper were presented at the Conference ``Poisson
2000'' held at CIRM, Luminy, France, June 26--30, 2000. I am grateful to a
number of people for helpful discussions and comments on this work:
especially, to A.Bolsinov, J.-P.Dufour, R.Flores Espinoza, J.Huebschmann,
M.V.Karasev, Y.Kosmann-Schwarzbach, P.Liber\-mann, K.Mackenzie, and A.Wein\-stein.

I would like also to thank G.D\'{a}vila Rasc\'{o}n for his help in the last
stages of preparation of this paper.

\section{Coupling tensors}

Let $\pi:\,E\rightarrow B$ be a fiber bundle (that is, a surjective
submersion). Let $\operatorname{Vert}=\ker\,(d\pi)\subset TE$ be the vertical
subbundle. Smooth sections of $\operatorname{Vert}$ form the Lie algebra of
vertical vector fields on the total space~$E$ which will be denoted by
$\mathcal{X}_{V}(E)$. Consider the annihilator $\operatorname{Vert}^{0}\subset
T^{\ast}E$ of the vertical subbundle. Sections of $\operatorname{Vert}^{0}$
are called horizontal $1$-forms. Denote by $\chi^{k}(E)=\Gamma(\Lambda^{k}TE)$
the space of $k$-vector fields on~$E$. A $k$-vector field $T\in\chi^{k}(E)$ is
said to be vertical if $\alpha\rfloor T=0$ for every horizontal $1$-form
$\alpha$$.$ The space of vertical $k$-vector fields will be denoted by
$\chi_{V}^{k}(E)$.

We say that a bivector field $\Pi\in\chi^{2}(E)$ is \textit{horizontally
nondegenerate }if for every $e\in E$ the antisymmetric bilinear form $\Pi
_{e}:\,T_{e}^{\ast}E\times T_{e}^{\ast}E\rightarrow\mathbb{R}$ is
nondegenerate on the subspace $\operatorname{Vert}_{e}^{0}\subset T_{e}^{\ast
}E$. In other words,
\begin{align}
\Pi^{\#}(\operatorname{Vert}^{0})\cap\operatorname{Vert} &  =\{0\},\\
\operatorname{rank}\Pi^{\#}(\operatorname{Vert}^{0}) &  =\dim B.
\end{align}
Here $\Pi^{\#}:\,T^{\ast}E\rightarrow TE$ is the vector bundle morphism
associated with~$\Pi$, $\Pi^{\#}(\alpha):=\alpha\rfloor\Pi$ ($\alpha\in
\Omega^{1}(E)$).

Recall that a bivector field $\Pi$ on $E$ is said to be a \textit{Poisson
tensor} if $\Pi$ satisfies the Jacobi identity \cite{Li}
\begin{equation}
\lbrack\![\Pi,\Pi]\!]_{E}=0.
\end{equation}
Here $[\![\cdot,\cdot]\!]_{E}$ denotes the \textit{Schouten bracket }for
multivector fields on the total space $E$. The corresponding Poisson bracket
is given by
\[
\{F,G\}_{\Pi}=\Pi(dF,dG)=\left\langle dG,\Pi^{\#}(dF)\right\rangle \text{.}%
\]
The correspondence $C^{\infty}(E)\ni F\mapsto\Pi^{\#}(dF)\in\mathcal{X(}E)$ is
a homorphism from the Poisson algebra $(C^{\infty}(E)$,$\{\,,\,\}_{\Pi})$ onto
the Lie algebra $\mathcal{X}^{Ham}\mathcal{(}E)$\ of Hamiltonian vector fields.

Our goal is to describe horizontally nondegenerate Poisson tensors on~$E$.
First, we formulate some preliminary facts.

\subsection{Geometric data}

Suppose we are given a horizontally nondegenerate bivector field $\Pi\in
\chi^{2}(E)$. By conditions (2.1), (2.2), we deduce that~$\Pi$ induces an
\textit{intrinsic Ehresmann connection} $\Gamma\in\Omega^{1}(E)\otimes
\mathcal{X}_{V}(E)$ whose horizontal subbundle $\operatorname{Hor}=\ker
\Gamma\subset TE$ is defined as the image of $\operatorname{Vert}^{0}$ under
the bundle map $\Pi^{\#}$,
\begin{equation}
\operatorname{Hor}:=\Pi^{\#}(\operatorname{Vert}^{0}).
\end{equation}
So, we have the splitting
\begin{equation}
TE=\operatorname{Hor}\oplus\operatorname{Vert}.
\end{equation}
Then $\mathbb{H=}\operatorname{id}-$ $\Gamma$ is the horizontal projection.
Let $\operatorname{Hor}^{0}\subset T^{\ast}E$ be the annihilator of the
horizontal subbundle. Sections of $\operatorname{Hor}^{0}$ are called vertical
$1$-forms. The set of $k$-vector fields $T\in\chi^{k}(E)$ such that
$\beta\rfloor T=0$ for every vertical $1$-form $\beta$ form the space of
horizontal $k$-vector fields denoted by $\chi_{H}^{k}(E)$. In particular,
$\mathcal{X}_{H}(E)=\chi_{H}^{1}(E)$ will denote the space of horizontal
vector fields on~$E$.

The splitting (2.5) induces a $C^{\infty}(B)$ homomorphism
\[
\operatorname{hor}:\,\mathcal{X}(B)\rightarrow\mathcal{X}_{H}(E)
\]
sending a smooth vector field $u$ on $B$ \ to a smooth section
$\operatorname{hor}(u) $ of $\ \operatorname{Hor}$ and satisfying
$L_{\operatorname{hor}(u)}(\pi^{\ast}f)=\pi^{\ast}(L_{u}f)$ for every $f\in
C^{\infty}(B)$. The vector field $\operatorname{hor}(u)$ is called the
\textit{horizontal lift} of a base vector field $u$, associated with the
connection~$\Gamma$. Notice that the flow $\operatorname{Fl}_{t} $ of
$\operatorname{hor}(u)$ is $\pi$-related with the flow $\varphi_{t}$ of $u$,
$\ \pi\circ\operatorname{Fl}_{t}=\varphi_{t}\circ\pi$.

It follows form (2.4) that subbundles $\operatorname{Hor}^{0}$ and
$\operatorname{Vert}^{0}$ are $\Pi$-orthogonal. Thus there is a unique
decomposition of~$\Pi$ into horizontal and vertical parts:
\begin{equation}
\Pi=\Pi_{H}+\Pi_{V},\qquad\text{where}\quad\Pi_{H}\in\chi_{H}^{2}(E),\quad
\Pi_{V}\in\chi_{V}^{2}(E).
\end{equation}
Consider the horizontal part $\Pi_{H}$. The horizontal nondegeneracy of $\Pi$
implies that the restriction
\begin{equation}
\overset{\circ}{\Pi}_{H}^{\#}:=\Pi_{H}^{\#}\big|_{\operatorname{Vert}^{0}%
}\,:\,\operatorname{Vert}^{0}\rightarrow\operatorname{Hor}%
\end{equation}
is a vector bundle isomorphism. Note that $\Pi_{H}$, as a horizontal bivector
field, is uniquely determined by $\overset{\circ}{\Pi}_{H}^{\#}$.

Consider the tensor product $\Omega^{k}(B)\otimes C^{\infty}(E)$ of
$C^{\infty}(B)$-modules. One can think of the elements of this space as
$k$-forms on the base~$B$ with values in the space $C^{\infty}(E)$, that is,
antisymmetric $k$-linear over $C^{\infty}(B)$ mappings $\mathcal{X}%
(B)\times\dots\times\mathcal{X}(B)\rightarrow C^{\infty}(E)$. Hence if
$\mathcal{F\in}\Omega^{k}(B)\otimes C^{\infty}(E)$ and $u_{1},...,u_{k}%
\in\mathcal{X}(B)$, then the restriction of $\mathcal{F(}u_{1},...,u_{k})$ to
the fiber\ $E_{b\text{ }}$depends only on $u_{1}(b),...,u_{k}(b)$ and we have
a $k$-linear over $\mathbb{R}$ mapping
\[
\mathcal{F}_{b}:T_{b}B\times...\times T_{b}B\rightarrow C^{\infty}(E_{b\text{
}})\text{.}%
\]
Moreover, there is a natural identification \ of $\Omega^{k}(B)\otimes
C^{\infty}(E)$ with the space $\Omega_{H}^{k}(E)$ of horizontal $k$-forms on
$E$: \ $\Omega^{k}(B)\otimes C^{\infty}(E)\ni{\mathcal{F}}\mapsto{\mathcal{F}%
}^{h}\in\Omega_{H}^{k}(E)$, where ${\mathcal{F}}^{h}$ is a horizontal $k$-form
defined \ by
\[
{\mathcal{F}}^{h}(Y_{1},...,Y_{k}):=\mathcal{F}_{\pi(e)}(d_{e}\pi
Y_{1},...,d_{e}\pi Y_{k})(e)
\]
for $Y_{1},...,Y_{k}\in T_{e}E$, $e\in E$. In particular, we have
${\mathcal{F}}^{h}(\operatorname{hor}(u_{1}),\dots,\operatorname{hor}%
(u_{k}))={\mathcal{F}}(u_{1},\dots,u_{k})$ and $(\omega\otimes1)^{h}=\pi
^{\ast}\omega$ for $\omega\in\Omega^{k}(B)$. We will say that $\mathcal{F\in
}\Omega^{2}(B)\otimes C^{\infty}(E)$ is \textit{nondegenerate} at a point
$e\in E$ if the restriction of ${\mathcal{F}}^{h}$ to the horisontal space
$\operatorname{Hor}_{e}\approx T_{e}E/\operatorname{Vert}_{e\text{ }}$ is a
nondegenerate bilinear form.

Let us associate to~$\Pi$ the $2$-form $\mathbb{F}\in\Omega^{2}(B)\otimes
C^{\infty}(E)$ defined by
\begin{equation}
\mathbb{F}(u_{1},u_{2}):=-\big\langle\big(\overset{\circ}{\Pi}_{H}^{\#}\big
)^{-1}\operatorname{hor}(u_{1}),\operatorname{hor}(u_{2})\big\rangle
\end{equation}
for $u_{1},u_{2}\in\mathcal{X}(B)$. Note that the $2$-form $\mathbb{F}$ is \textit{nondegenerate,}%

\begin{equation}
u\rfloor\mathbb{F}=0\qquad\text{for $u\in\mathcal{X}(B)$ \textit{implies}
$u=0$}.
\end{equation}
Here the interior product $u\rfloor\mathbb{F}$ is an element of the space
$\Omega^{1}(B)\otimes C^{\infty}(E)$.

Now we claim that horizontally nondegenerate bivector fields on~$E$ can be
parametrized by some geometric data. By \textit{geometric data} we mean a
triple $(\Gamma,\mathcal{V},\mathbb{F})$ consisting of

\begin{itemize}
\item  an Ehresmann connection~$\Gamma$ on~$\pi$;

\item  a vertical bivector field $\mathcal{V}\in\chi^{2}_{V}(E)$;

\item  a nondegenerate $C^{\infty}(E)$-valued $2$-form $\mathbb{F}\in
\Omega^{2}(B)\otimes C^{\infty}(E)$ on the base~$B$.
\end{itemize}

\textit{Direct mapping $\Pi\mapsto(\Gamma,\mathcal{V},\mathbb{F})$}. We
associate to a given horizontally nondegenerate bivector field~$\Pi$ on~$E$
the geometric data $(\Gamma,\mathcal{V},\mathbb{F})$, where

\begin{itemize}
\item $\Gamma:TE\rightarrow\operatorname{Vert}$ is the projection along \ the
subbundle (2.4);

\item $\mathcal{V}=\Pi_{V}$ is the vertical part of $\Pi$ in (2.6);

\item $\mathbb{F}$ is the $2$-form in (2.8).
\end{itemize}

\textit{Inverse mapping $(\Gamma,\mathcal{V},\mathbb{F})\mapsto\Pi$}. Taking
geometric data $(\Gamma,\mathcal{V},\mathbb{F})$, we introduce the
horizontally nondegenerate bivector field
\begin{equation}
\Pi=\tau_{\Gamma}(\mathbb{F})+\mathcal{V},
\end{equation}
where the $\Gamma$-dependent correspondence $\mathbb{F}\mapsto\tau_{\Gamma
}(\mathbb{F})\in\chi_{H}^{2}(E)$ is defined in the following way. The
nondegenerate $2$-form $\mathbb{F}\in\Omega^{2}(B)\otimes C^{\infty}(E)$
induces the vector bundle isomorphism
\begin{equation}
\mathbb{F}^{\,\flat}:\,\operatorname{Hor}\rightarrow\operatorname{Vert}%
^{0}\approx\operatorname{Hor}^{\ast}%
\end{equation}
such that
\begin{equation}
\langle\mathbb{F}^{\,\flat}(\operatorname{hor}(u_{1})),\operatorname{hor}%
(u_{2})\rangle=\mathbb{F}(u_{1},u_{2})
\end{equation}
for every $u_{1},u_{2}\in\mathcal{X}(B)$. Then the horizontal bivector field
$\tau_{\Gamma}(\mathbb{F})$ is determined by the condition
\begin{equation}
\tau_{\Gamma}(\mathbb{F})(\beta_{1},\beta_{2})=-\langle\beta_{2}%
,(\mathbb{F}^{\,\flat})^{-1}\beta_{1}\rangle
\end{equation}
for all $\beta_{1},\beta_{2}\in\Gamma(\operatorname{Vert}^{0})$.

\subsection{Coupling tensors}

Now we can try to rewrite the Jacobi identity (2.3) for a horizontally
nondegenerate bivector field in terms of its geometric data. To formulate the
result, we recall some definitions.

If $\Gamma$ is an Ehresmann connection on a fiber bundle $\pi:\, E\to B$, then
the curvature form is a vector valued $2$-form $\operatorname{Curv}^{\Gamma
}\in\Omega^{2}(B,\operatorname{Vert})\approx\Omega^{2}(B)\otimes
\mathcal{X}_{V}(E)$ on the base defined as
\[
\operatorname{Curv}^{\Gamma}(u_{1},u_{2}):=-\Gamma\big([\operatorname{hor}%
(u_{1}),\operatorname{hor}(u_{2})]\big)%
\]
for $u_{1},u_{2}\in\mathcal{X}(B)$.

The connection $\Gamma$ induces the \textit{covariant exterior derivative}
\cite{GHV}
\[
\partial_{\Gamma}:\Omega^{k}(B)\otimes C^{\infty}(E)\rightarrow\Omega
^{k+1}(B)\otimes C^{\infty}(E)
\]
taking a $k$-form $\mathcal{F}$ to a $(k+1)$-form $\partial_{\Gamma
}\mathcal{F}$, which at vector fields $u_{0},u_{1},\dots,u_{k}\in
\mathcal{X}(B)$ is:
\begin{align}
&  (\partial_{\Gamma}\mathcal{F})(u_{0},u_{1},\dots,u_{k}):=\sum_{i=0}%
^{k}(-1)^{i}L_{\operatorname{hor}(u_{i})}\mathcal{F}(u_{0},u_{1},\dots,\hat
{u}_{i},\dots,u_{k})\nonumber\\
&  \qquad+\sum_{0\leq i<j\leq k}(-1)^{i+j}\mathcal{F}\big([u_{i},u_{j}%
],u_{0},u_{1},\dots,\hat{u}_{i},\dots,\hat{u}_{j},\dots,u_{k}\big).
\end{align}
Notice that $\partial_{\Gamma}$ is a coboundary operator, $\partial_{\Gamma
}^{2}=0$ if and only if $\operatorname{Curv}^{\Gamma}=0$, that is, $\Gamma$ is
\textit{flat}. Moreover, we have $(\partial_{\Gamma}{\mathcal{F}}%
)^{h}=\mathbb{H}_{\ast}\circ d({\mathcal{F}}^{h})$. Here $\mathbb{H}_{\ast
}:\Omega^{k}(E)\rightarrow\Omega_{H}^{k}(E)$ is the horizontal projection and
$d$ is the usual differential of forms. In particular,
\begin{equation}
\partial_{\Gamma}(\omega\otimes1)=d\omega\otimes1
\end{equation}
for every $k$-form $\omega$ on the base $B$.

\begin{thm}
Let $\pi:\,E\to B$ be a fiber bundle. A horizontally nondegenerate bivector
field $\Pi\in\chi^{2}(E)$ is a Poisson tensor if and only if its geometric
data $(\Gamma,\mathcal{V},\mathbb{F})$ satisfy the following conditions:

\emph{(i)} the vertical part $\mathcal{V}$ of $\Pi$ is a Poisson tensor,
\begin{equation}
[\![\mathcal{V},\mathcal{V}]\!]_{E}=0;
\end{equation}

\emph{(ii)} the connection $\Gamma$ preserves $\mathcal{V}$, that is, for
every $u\in\mathcal{X}(B)$ the horizontal lift $\operatorname{hor}(u)$ is an
infinitesimal automorphism of $\mathcal{V}$,
\begin{equation}
L_{\operatorname{hor}(u)}\mathcal{V}\equiv[\![\operatorname{hor}%
(u),\mathcal{V}]\!]_{E}=0;
\end{equation}

\emph{(iii)} the nondegenerate $2$-form $\mathbb{F}\in\Omega^{2}(B)\otimes
C^{\infty}(E)$ satisfies
\begin{equation}
\partial_{\Gamma}\mathbb{F}=0,
\end{equation}
and the ``curvature identity''
\begin{equation}
\operatorname{Curv}^{\Gamma}(u_{1},u_{2})=\mathcal{V}^{\#}d\mathbb{F}%
(u_{1},u_{2})
\end{equation}
for $u_{1},u_{2}\in\mathcal{X}(B)$.
\end{thm}

The proof is a direct verification under the use of the Poisson--Ehresmann
calculus. For a symplectic version of Theorem 2.1. see \cite{GLS}.

So, if geometric data $(\Gamma,\mathcal{V},\mathbb{F})$ satisfy conditions
(2.16)--(2.19), then formula (2.10) determines a Poisson tensor~$\Pi$ on~$E$
which will be called a \textbf{coupling tensor} associated with $(\Gamma
,\mathcal{V},\mathbb{F})$.

The hypotheses in Theorem~2.1 have the following interpretations. Suppose we
are given some geometric data $(\Gamma,{\mathcal{V}},\mathbb{F})$ satisfying
conditions (2.16)--(2.19).

It follows from (2.16) that $(E,\mathcal{V)}$ is a Poisson manifold with
vertical Poisson tensor $\mathcal{V}$. Let $\operatorname{Casim}%
_{{\mathcal{V}}}(E)$ be the space of \textit{Casimir functions} on
$(E,\mathcal{V)}$, that is, the center of the Poisson algebra $(C^{\infty
}(E),\{\,,\,\}_{\mathcal{V}})$. Clearly, $\pi^{\ast}C^{\infty}(B)\subset
\operatorname{Casim}_{{\mathcal{V}}}(E).$ Every fiber $E_{b}=\pi^{-1}(b)$
($b\in B$) is a \textit{Poisson submanifold} of $(E,{\mathcal{V}})$ which
carries a unique Poisson structure ${\mathcal{V}}_{b}^{\text{fib}}$ with
property: the inclusion $E_{b}\hookrightarrow E$ is a Poisson mapping (see
\cite{We4}). Thus, ${\mathcal{V}}$ induces a smooth field of Poisson
structures on the fibers: $B\ni b\mapsto{\mathcal{V}}_{b}^{\text{fib}}\in
\chi^{2}(E_{b})$ called a \textit{fiberwise Poisson structure}. Notice that if
we start with a locally trivial Poisson fiber bundle (the typical fiber is a
Poissin manifold), then the fiberwise Poisson structure induces a unique
compatible vertical Poisson tensor.

Condition (2.17) means that the horizontal lift $\operatorname{hor}(u)$ of
every base vector field~$u$ is a \textit{Poisson vector field} (an
\textit{infinitesimal Poisson automorphism}) of the vertical Poisson
structure~${\mathcal{V}}$. The Ehresmann connection $\Gamma$ is compatible
with the fiberwise Poisson structure in the sense that the (local) flow
$\operatorname{Fl}_{t} $ of $\operatorname{hor}(u)$ is a fiber preserving
Poisson morphism. In other words, the parallel transport associated with the
connection $\Gamma$ preserves the fiberwise Poisson structure. Such a
connection is called a \textit{Poisson connection} on a bundle of Poisson manifolds.

Denote by ${\mathcal{X}}_{V}^{\text{Poiss}}(E)$ the Lie algebra of vertical
Poisson vector fields and by ${\mathcal{X}}_{V}^{\text{Ham}}(E)$ the Lie
subalgebra of vertical Hamiltonian vector fields on $(E,{\mathcal{V}})$. Then
${\mathcal{X}}_{V}^{\text{Ham}}(E)$ is an ideal in ${\mathcal{X}}%
_{V}^{\text{Poiss}}(E)$. Consider the quotient space
\begin{equation}
{\mathcal{H}}_{V}^{1}(E;{\mathcal{V}})={\mathcal{X}}_{V}^{\text{Poiss}%
}(E)/{\mathcal{X}}_{V}^{\text{Ham}}(E)\text{.}%
\end{equation}
Notice that the symplectic leaf of ${\mathcal{V}}$ through a point $e\in E$
coincides with the symplectic leaf of the Poisson structure\ ${\mathcal{V}%
}_{b}^{\text{fib}}$ on the fiber $E_{\pi(e)}$. Hence every Hamiltonian vector
field on $(E,{\mathcal{V}})$ is an element of ${\mathcal{X}}_{V}^{\text{Ham}%
}(E)$. Moreover, every Poisson vector field of ${\mathcal{V}}$ is represented
as a sum of a horozontal lift of some base field and an element of
${\mathcal{X}}_{V}^{\text{Poiss}}(E)$. From here we deduce: the first Poisson
cohomology space\cite{Li,KM,Va} of ${\mathcal{V}}$ is isomorphic to the direct
sum $\mathcal{X}(B)\oplus{\mathcal{H}}_{V}^{1}(E;{\mathcal{V}})$.

One can show that conditions (2.17) implies the property: for every
$u_{1},u_{2}\in{\mathcal{X}}(B)$ the curvature vector field
$\operatorname{Curv}^{\Gamma}(u_{1},u_{2})$ is a vertical Poisson vector
field,
\begin{equation}
\operatorname{Curv}^{\Gamma}(u_{1},u_{2})\in{\mathcal{X}}_{V}^{\text{Poiss}%
}(E).
\end{equation}
The curvature identity (2.19) leads to the stronger requirement:
$\operatorname{Curv}^{\Gamma}(u_{1},u_{2})$ is a vertical Hamiltonian vector
field with the Hamiltonian function $\mathbb{F}(u_{1},u_{2})$,
\begin{equation}
\operatorname{Curv}^{\Gamma}(u_{1},u_{2})\in{\mathcal{X}}_{V}^{\text{Ham}}(E)
\end{equation}
and hence the equivalence class of $\operatorname{Curv}^{\Gamma}(u_{1},u_{2})$
in ${\mathcal{H}}_{V}^{1}(E,{\mathcal{V}})$ is trivial.

Remark also that conditions (2.18) and (2.19) are independent in general.
Indeed, the curvature identity (2.19) and the Bianchi identity for the
curvature form of $\Gamma$ imply only that
\begin{equation}
\partial_{\Gamma}\mathbb{F}\in\Omega^{2}(B)\otimes\operatorname{Casim}%
_{{\mathcal{V}}}(E).
\end{equation}
\ Consider the following two ''extreme'' cases.

\begin{exam}
[Flat Poisson bundles]\textrm{Suppose we start with a Poisson bundle over
symplectic base: $(E\xrightarrow{\pi}B,{\mathcal{V}},\omega)$, where
${\mathcal{V}}$ is a vertical Poisson tensor and $\omega$ is a base symplectic
structure. Then we can assign to every \textit{flat} Poisson connection
$\Gamma$ on $\pi$ the coupling tensor $\Pi^{\Gamma}$ defined by (2.10), where
$\mathbb{F}=\omega\otimes1$. In this case, (2.19) holds because of the
flatness, $\operatorname{Curv}^{\Gamma}=0$. And (2.18) follows from the
closedness of $\omega$. The horizontal part $\Pi_{H}^{\Gamma}$ is just lifting
of the nondegenerate Poisson structure on $(B,\omega)$ via $\Gamma$. In the
nonflat case, to satisfy the curvature identity (2.19) we have to deform the
symplectic structure on the base. }
\end{exam}

\begin{exam}
[Coupling forms \cite{St,GLS}]\textrm{Under the same starting\break point as
in Example 2.1, assume also that the Poisson structure ${\mathcal{V}}_{b}$ is
\textit{nondegenerate} on each fiber $E_{b}$. Then ${\mathcal{V}}$ induces a
fiberwise symplectic structure }$B\ni b\mapsto\sigma_{b}\in\Omega^{2}(E_{b}%
)$\textrm{\ and the bundle $E$ becomes a \textit{symplectic fiber}%
\textbf{\ }\textit{bundle} $(E,\sigma)$ over a symplectic base (}$B,\omega
)$\textrm{. In this case, every Poisson connection $\Gamma$ on $E$ is also
symplectic, that is, the parallel transport preserves the fiberwise symplectic
structure $\sigma$. Furthermore, }$\operatorname{Casim}_{{\mathcal{V}}%
}(E)\approx C^{\infty}(B).$\textrm{\ Let us make the extra assumption: the
fiber bundle $\pi:\,E\rightarrow B$ } is locally trivial and \textrm{the
typical fiber is \textit{compact, connected} and \textit{simply connected}.
Then we have (see \cite{GLS})
\begin{equation}
{\mathcal{H}}_{V}^{1}(E;{\mathcal{V}})=0.
\end{equation}
Let $\Gamma$ be a Poisson connection. It follows from (2.21) and (2.24) that
(2.22) holds. The problem now is to find $\mathbb{F}$ in (2.19) satisfying
also condition (2.18). Taking into account our assumption, introduce a
}$C^{\infty}(B)$ linear\textrm{\ mapping
\[
{\mathcal{X}}_{V}^{\text{Ham}}(E)\ni Z\mapsto{\mathbf{m}}(Z)\in C^{\infty}(E)
\]
which is determined by conditions
\begin{align}
Z\rfloor\sigma_{b} &  =-d{\mathbf{m}}_{b}(Z)\qquad\text{on}\quad E_{b},\\
\int_{E_{b}}{\mathbf{m}}_{b}(Z)\sigma_{b}^{n} &  =0,\qquad2n=\dim
(\text{\textit{fiber}})
\end{align}
for every $b\in B$. Here ${\mathbf{m}}_{b}(Z)={\mathbf{m}}(Z)\big|_{E_{b}}$
and $\sigma_{b}^{n}=\frac{1}{n!}\sigma_{b}\wedge\dots\wedge\sigma_{b}$
($n$-times) is the volume form. Then we claim that the formula
\[
\mathbb{F}^{\Gamma}(u_{1},u_{2})=\pi^{\ast}\omega(u_{1},u_{2})+{\mathbf{m}%
}(\operatorname{Curv}^{\Gamma}(u_{1},u_{2}))
\]
defines just the desired form $\mathbb{F}^{\Gamma}\in\Omega^{2}(B)\otimes
C^{\infty}(E)$ satisfying (2.18) and (2.19). Indeed, (2.18) holds
automatically. The symplecticity of $\Gamma$ and the normalization condition
(2.26) imply that the ``collective Hamiltonian'' ${\mathbf{m}}$ is $\Gamma
$-invariant, \textbf{m}$((\operatorname{Fl}_{t}^{-1})^{\ast}%
Z)=\operatorname{Fl}_{t}^{\ast}(\mathbf{m}(Z))$, where $\operatorname{Fl}_{t}$
is the flow of $\operatorname{hor}(u)$. This leads to (2.19). Denote by
$\sigma^{v}$ the vertical $2$-form on $E$ which coincides with $\sigma_{b}$ on
each fiber $E_{b}$. Then the closed $2$-form
\begin{equation}
\Omega^{\Gamma}=(\mathbb{F}^{\Gamma})^{h}+\sigma^{v}%
\end{equation}
is called a \textbf{coupling form} associated with the symplectic connection
$\Gamma$ \cite{GLS}. The closedness of $\Omega^{\Gamma}$ is just equivalent to
conditions (2.18) and (2.19) for $\mathbb{F}=\mathbb{F}^{\Gamma}$. In a domain
where $\mathbb{F}^{\Gamma}$ is nondegenerate, $\Omega^{\Gamma}$ is symplectic
and its nondegenerate Poisson structure is the coupling tensor $\Pi^{\Gamma}$
generated by the triple $(\Gamma,{\mathcal{V}},\mathbb{F}^{\Gamma})$. }
\end{exam}

\section{Neighborhood equivalence}

Here we show that coupling tensors naturally appear in the classification
problem of Poisson structures near a single symplectic leaf.

\subsection{Geometric splitting}

Let $\nu:\,{\mathcal{N}}\to B$ be a fiber bundle over a connected base $B$.
Suppose we have a cross-section ${{\mathbf{s}}}:\,B\to{\mathcal{N}}$ of $\nu$.

We say that a Poisson tensor $\Pi$ on ${\mathcal{N}}$ is compatible with
section ${\mathbf{s}}$, or shortly, ${\mathbf{s}}$-\textbf{compatible} if
${\mathbf{s}}(B)\text{ is a \textit{symplectic leaf\/} of }\Pi$.

\begin{prop}
Let $\Pi\in\chi^{2}({\mathcal{N}})$ be a ${\mathbf{s}}$-compatible Poisson
tensor. Then there exists a tubular neighborhood $E$ of ${\mathbf{s}}(B)$ in
${\mathcal{N}}$ such that $\Pi$ is a coupling tensor on $E$. In particular,
there is an \textit{intrinsic Ehresmann connection} $\Gamma$ on $\pi=\nu\big
|_{E}:\,E\rightarrow B$ which induces a unique decomposition
\begin{equation}
\Pi=\Pi_{H}+\Pi_{V}\quad\text{on}\quad E,
\end{equation}
where

\begin{itemize}
\item $\Pi_{H}\in\chi^{2}_{H}(E)$ is a $\Gamma$-horizontal bivector field,

\item $\Pi_{V}\in\chi^{2}_{V}(E)$ is a vertical Poisson tensor such that
\begin{equation}
\operatorname{rank}\Pi_{V}=0\quad\text{at every point in}\quad{\mathbf{s}}(B).
\end{equation}
\end{itemize}

The connection $\Gamma$ is determined by \emph{(2.4)}.
\end{prop}

\begin{proof}
Since ${\mathbf{s}}(B)$ is a symplectic leaf of $\Pi$, the bivector field
$\Pi$ is nondegenerate on subspaces $\operatorname{Vert}_{e}^{0}\subset
T_{e}^{\ast}N$ at points $e\in{\mathcal{N}}$ sufficiently close to
${\mathbf{s}}(B)$. Hence $\Pi$ is a horizontally nondegenerate on a tubular
neighborhood $E$ of ${\mathbf{s}}(B)$ in ${\mathcal{N}}$. By Theorem~2.1,
$\Pi$ is a coupling tensor on $E$ associated with geometric data
$(\Gamma,{\mathcal{V}}=\Pi_{V},\mathbb{F})$. Here the intrinsic connection
$\Gamma$ and the bivector fields $\Pi_{H}$, ${\mathcal{V}}$ are defined by
(2.4) and (2.6) respectively. Property (3.2) follows form ${\mathbf{s}}%
$-compatibility assumption.
\end{proof}

\begin{rem}
\textrm{Each fiber $E_{b}$ over $b\in B$ inherits from }$\Pi$ \textrm{a
Poisson structure in a neighborhood of ${\mathbf{s}}(b)$ called the
\textit{transverse Poisson structure} at the point }$b$ \textrm{\cite{We4}. By
the splitting theorem, transverse Poisson structure is independent of the
choice of a point ${\mathbf{s}}(b)$ up to local isomorphism. The vertical part
$\Pi_{V}$ in (3.1) gives rise to a fiberwise Poisson structure which fit
together local transverse Poisson structures. }
\end{rem}

Note that the maximal domain, where splitting (3.1) holds, consists of the
points $e\in{\mathcal{N}}$ such that $\operatorname{rank}_{e}\Pi
^{\#}(\operatorname{Vert}^{0})=\dim B$. A given ${\mathbf{s}}$-compatible
Poisson structure $\Pi$ with geometric data $(\Gamma,\Pi_{V},\mathbb{F})$
possesses the following properties on $E$.

(i) The symplectic structure $\omega$ on $\mathbf{s}(B)$ is
\begin{equation}
\omega=\mathbb{F}^{h}\big|_{{\mathbf{s}}(B)}.
\end{equation}
The horizontal distribution $\operatorname{Hor}$ associated with $\Gamma$ is
tangent to ${\mathbf{s}}(B)$,
\begin{equation}
T_{e}{\mathbf{s}}(B)=\operatorname{Hor}_{e}\qquad\text{for}\quad
e\in{\mathbf{s}}(B).
\end{equation}
and hence
\begin{equation}
\operatorname{Curv}^{\Gamma}(u_{1},u_{2})=0\quad\text{on}\quad{\mathbf{s}}(B).
\end{equation}
The projection $\pi:E\rightarrow B$ is a Poisson morphism if and only if the
curvature of $\Gamma$ is zero, $\operatorname{Curv}^{\Gamma}=0$ on $E$. In the
flat case, the symplectic leaf \ $\mathbf{s}(B)$ is an integral leaf of the
integrable horizontal distribution $\operatorname{Hor}\subset TE$. Hence,
there is the holonomy of $\mathbf{s}(B)$ (as a leaf of the corresponding
foliation), called the \textit{strict Poisson holonomy} of the leaf \cite{Fe}.

(ii) For every $f\in C^{\infty}(B)$ the Hamiltonian vector field of the pull
back $\pi^{\ast}f$ is horizontal,
\begin{equation}
\Pi^{\#}(\pi^{\ast}df)=\sum_{i}\Pi(\pi^{\ast}df,\pi^{\ast}d\xi^{i}%
)\operatorname{hor}(\partial/\partial\xi^{i}),
\end{equation}
where $(\xi^{i})$ are (local) coordinates on the base $B$. \ Let
$\mathcal{D}=\operatorname{Ann}(\ker\Pi^{\#})$ be the \textit{characteristic
distribution} of $\Pi$ on $E$. Let ${\mathcal{V}}_{b}^{\text{fib}}\in\chi
^{2}(E_{b})$ be the Poisson tensor on the fiber $E_{b}$ generated by $\Pi_{V}%
$, and let ${\mathcal{D}}_{b}^{\text{fib}}$ be the characteristic distribution
of ${\mathcal{V}}_{b}^{\text{fib}}$. Then for every $e\in E$ we
have$\ {\mathcal{D}}_{_{e}}=\operatorname{Hor}_{e}\oplus({\mathcal{D}}%
_{b}^{\text{fib}})_{e}$. Hence the rank of the Poisson structure $\Pi$ at $e$
is $\operatorname{rank}_{e}\Pi=\dim B+\operatorname{rank}_{e}{\mathcal{V}}%
_{b}^{\text{fib}},$ $\ b=\pi(e).$ If $(\mathcal{S},\Omega)$ is a symplectic
leaf of $(E,\Pi)$ with symplectic structure $\Omega$, then at every $e\in E$
we have the decomposition
\[
\Omega_{e}=(\mathbb{F}^{h})_{e}\oplus\sigma_{e}\text{,}%
\]
where $\sigma$ ia the symplectic form on the leaf of $(E_{\pi(e)\text{,}%
}{\mathcal{V}}_{\pi(e)}^{\text{fib}})$ passing through the point $e$. The
Poisson tensor $\Pi$ is of constant rank on $E$ if and only if $\Pi_{V}%
\equiv0$. In this case, $\Gamma$ is flat.

(iii) Let $f:\,\widetilde{E}\rightarrow E$ be a fiber preserving
diffeomorphism from an open neighborhood $\widetilde{E}$ of ${\mathbf{s}}(B)$
onto $E$ that descends the identity map on ${\mathbf{s}}(B)$. Then $f^{\ast
}\Pi$ is a ${\mathbf{s}}$-compatible Poisson tensor on $\widetilde{E}$ with
the intrinsic connection $f^{\ast}\Gamma$ and the splitting \ $f^{\ast}%
\Pi=(f^{\ast}\Pi)_{H}+(f^{\ast}\Pi)_{V}=f^{\ast}(\Pi_{H})+f^{\ast}(\Pi_{V}).$

As a direct consequence of Proposition 3.1, we get the fact: in a tubular
neighborhood of a closed symplectic leaf every Poisson structure is realized
as a coupling tensor (see Section~5). Thus, the problem on the neighborhood
equivalence between Poisson structures near a common symplectic leaf is
reduced to the investigation of coupling tensors over a compatible cross-section.

\subsection{Neighborhood equivalence}

Let $\pi:E\rightarrow B$ be a fiber bundle over a connected base. Suppose we
have two coupling tensors $\Pi$ and $\widetilde{\Pi}$ on $E$ associated with
geometric data $(\Gamma,{\mathcal{V}},\mathbb{F})$ and $(\widetilde{\Gamma
},\widetilde{{\mathcal{V}}},\widetilde{\mathbb{F}})$, respectively. Assume
that $\Pi$ and $\widetilde{\Pi}$ are compatible with a cross-section
${\mathbf{s}}:\,B\rightarrow E$ and
\begin{equation}
\mathbb{F}(u_{1},u_{2})\big|_{{\mathbf{s}}(B)}=\widetilde{\mathbb{F}}%
(u_{1},u_{2})\big|_{{\mathbf{s}}(B)}%
\end{equation}
for $u_{1},u_{2}\in{\mathcal{X}}(B)$. Condition (3.7) means that the
symplectic structures on ${\mathbf{s}}(B)$ with respect to Poisson structures
$\Pi$ and $\widetilde{\Pi}$ coincide.

We say that the geometric data $(\Gamma,{\mathcal{V}},\Lambda)$ and
$(\widetilde{\Gamma},\widetilde{{\mathcal{V}}},\widetilde{\Lambda})$ are
\textbf{equivalent} over ${\mathbf{s}}(B)$ if there exist open neighborhoods
$\mathcal{E}$ and $\widetilde{\mathcal{E}}$ of ${\mathbf{s}}(B)$ in $E$ and a
pair $(g,\phi)$ consisting of

\begin{itemize}
\item  a fiber preserving diffeomorphism $g:\,\mathcal{E}\rightarrow
\widetilde{\mathcal{E}}$ ( $\pi\circ g=\pi$) \ such that $g\circ
\mathbf{s}=\mathbf{s}$;

\item  a base $1$-form $\phi\in\Omega^{1}(B)\otimes C^{\infty}(\mathcal{E})$
\end{itemize}

which satisfy the relations:
\begin{align}
g^{\ast}\widetilde{{\mathcal{V}}} &  ={\mathcal{V}},\\
g^{\ast} &  =\Gamma-({\mathcal{V}}^{\#}\,d\phi)^{h},\\
g^{\ast}\widetilde{\mathbb{F}} &  =-\partial_{\Gamma}\phi-\frac{1}{2}%
\{\phi\wedge\phi\}_{{\mathcal{V}}}.
\end{align}
Here ${\mathcal{V}}^{\#}\,d\phi$ is an element of the space $\Omega
^{1}(B)\otimes\mathcal{X}_{V}^{Ham}(\mathcal{E})$ determined by $\mathcal{(}%
{\mathcal{V}}^{\#}\,d\phi)(u)={\mathcal{V}}^{\#}\,d\phi(u)$ and $\{\phi
_{1}\wedge\phi_{2}\}_{{\mathcal{V}}}$ denotes an element of $\Omega
^{2}(B)\otimes C^{\infty}(\mathcal{E})$ given by
\[
\frac{1}{2}\{\phi_{1}\wedge\phi_{2}\}_{{\mathcal{V}}}(u_{1},,u_{2}%
):={\mathcal{V}}\big(d\phi_{1}(u_{1}),d\phi_{2}(u_{2})\big)-{\mathcal{V}}\big
(d\phi_{1}(u_{2}),d\phi_{2}(u_{1})\big)%
\]
for $\phi_{1},\phi_{2}\in\Omega^{1}(B)\otimes C^{\infty}(\mathcal{E})$ and
$u_{1},u_{2}\in{\mathcal{X}}(B)$.

\begin{thm}
Let $\Pi$ and $\widetilde{\Pi}$ be two ${\mathbf{s}}$-compatible coupling
tensors satisfying condition \emph{(3.7)}. If the corresponding geometric data
$(\Gamma,{\mathcal{V}},\mathbb{F})$ and $(\widetilde{\Gamma},\widetilde
{{\mathcal{V}}},\widetilde{\mathbb{F}})$ are equivalent over ${\mathbf{s}}%
(B)$, then there exist neighborhoods ${\mathcal{O}}$ and $\widetilde
{{\mathcal{O}}}$ of ${\mathbf{s}}(B)$ in $E$ and a diffeomorphism
${\mathbf{f}}:\,{\mathcal{O}}\rightarrow\widetilde{{\mathcal{O}}}$ such that
\begin{align}
{\mathbf{f}}\circ{\mathbf{s}} &  ={\mathbf{s}},\\
{\mathbf{f}}^{\ast}\widetilde{\Pi} &  =\Pi.
\end{align}
\end{thm}

\begin{proof}
We will use a contravariant analog of the homotopy method due to Moser
\cite{Mos} and Weinstein \cite{We1} (see also \cite{GSt,LMr}.

\textbf{Step 1}. \textit{Homotopy between coupling tensors}. By the property
(iii) in section 3.1, without loss of generality we can assume that the
vertical parts of $\Pi$ and $\widetilde{\Pi}$ coincide, $\widetilde
{{\mathcal{V}}}={\mathcal{V}}$ $\ $on $\mathcal{E}$ and $g=\operatorname{id}$.
By the ${\mathbf{s}}$-compatibility assumption we deduce that
$\operatorname{rank}{\mathcal{V}}=0$ at ${\mathbf{s}}(B)$. It follows from
this property and (3.7) that we\ can choose $\phi$ in (3.9), (3.10) so that
\begin{equation}
\phi(u)\big|_{{\mathbf{s}}(B)}=0,\qquad\text{for all }u\in{\mathcal{X}}(B).
\end{equation}
Consider the following $t$-parameter families of forms:
\begin{align}
\Gamma_{t} &  =\Gamma-t({\mathcal{V}}^{\#}d\phi)^{h}\in\Omega^{1}%
(E)\otimes{\mathcal{X}}_{V}(\mathcal{E}),\\
\mathbb{F}_{t} &  =\mathbb{F}-t\partial_{\Gamma}\phi-\frac{t^{2}}{2}%
\{\phi\wedge\phi\}_{{\mathcal{V}}}\in\Omega^{2}(B)\otimes C^{\infty
}(\mathcal{E}).
\end{align}
Then $\Gamma_{t}$ is a time-dependent connection $1$-form on $\mathcal{E}$. By
(3.13) and the nondegeneracy of $\mathbb{F}$, there is a neighborhood
$\mathcal{E}_{0}$ of ${\mathbf{s}}(B)$ in $E$ such that $\mathbb{F}_{t}$ is
nondegenerate on $\mathcal{E}_{0}$ for all $t\in\lbrack0,1]$. This means that
for every $t\in\lbrack0,1]$ and $e\in$ $\mathcal{E}_{0}$ the horizontal lift
$($ $\mathbb{F}_{t})^{h}$ induces a nondegenerate bilinear form on the
qoutient \ space $T_{e}E/Vert_{e}$. Moreover, we observe that the triple
$(\Gamma_{t},{\mathcal{V}},\mathbb{F}_{t})$ defines a geometric data on
$\mathcal{E}_{0}$ satisfying conditions (2.16)--(2.19) for every $t\in
\lbrack0,1]$. Thus the time-dependent coupling tensor $\Pi_{t}$ ($t\in
\lbrack0,1]$) associated to $(\Gamma_{t},{\mathcal{V}},\mathbb{F}_{t}) $ gives
a homotopy from $\Pi$ to $\widetilde{\Pi}$, $\ \Pi_{t}\big|_{t=0}=\Pi,$
$\Pi_{t}\big|_{t=1}=\widetilde{\Pi}$.

\textbf{Step 2}. \textit{Homological equation}. By the nondegeneracy of
$\mathbb{F}_{t}$ on $\mathcal{E}_{0}$, there exists a unique solution
$X_{t}\in{\mathcal{X}}(B)\otimes C^{\infty}(\mathcal{E}_{0})$ of the following
equation
\begin{equation}
X_{t}\rfloor\mathbb{F}_{t}=\phi.
\end{equation}
Clearly
\begin{equation}
\langle df,X_{t}\rangle\big|_{{\mathbf{s}}(B)}=0\quad\text{for every}\quad
f\in C^{\infty}(B).
\end{equation}
One can associate to $X_{t}$ the time-dependent horizontal vector field
$X_{t}^{h}\in{{\mathcal{X}}}_{H}(\mathcal{E}_{0})$ defined by
\begin{equation}
X_{t}^{h}(\pi^{\ast}f)=\langle df,X_{t}\rangle
\end{equation}
for all $f\in C^{\infty}(B)$. Property (3.17) implies
\begin{equation}
X_{t}^{h}\big|_{{\mathbf{s}}(B)}=0.
\end{equation}

\begin{lem}
$X^{h}_{t}$ satisfies the equation
\begin{equation}
L_{X^{h}_{t}}\Pi_{t}+\frac{\partial}{\partial t}\Pi_{t}=0\qquad(t\in[0,1]).
\end{equation}
\end{lem}

Here $L_{X}\Pi$ is the Lie derivative of a bivector field $\Pi$ along a vector
field $X$, that is, the Schouten bracket $[\![X,\Pi]\!]_{E}$. The proof of
Lemma~3.1 is given in Appendix~A.

\textbf{Step 3}. Let $\Phi_{t}$ be the flow of the time-dependent horizontal
vector field $X_{t}^{h}$: $\frac{d}{dt}\Phi_{t}=X_{t}^{h}\circ\Phi_{t}$,
$\Phi_{0}=id$. By (3.20) and the usual properties of the Lie
derivative\ \ (see, for example, \cite{KM,LMr,Va} ) we get $\Phi_{t}^{\ast}%
\Pi_{t}=\Pi$. Because of (3.19) for every $e\in{\mathbf{s}}(B)$, we have
$\Phi_{t}(e)=e$ for all $t\in\lbrack0,1]$. Hence there exists a neighborhood
${\mathcal{O}}$ of ${\mathbf{s}}(B)$ in $E_{0}$ that lies in the domain of the
flow $\Phi_{t}$ for $t\in\lbrack0,1]$. Finally, the time 1 flow $\Phi_{1}$ of
$X_{t}^{h}$ generates a diffeomorphism ${\mathbf{f}}:\,{\mathcal{O}%
}\rightarrow\widetilde{{\mathcal{O}}}$ satisfying (3.11), (3.12).
\end{proof}

\begin{rem}
$\mathrm{\ If}$ ${\mathbf{s}}(B)$ is a regular symplectic leaf, then
$\mathbb{F}=\widetilde{\mathbb{F}}=0$ and ${\mathcal{V=}}0$. Condition (3.9)
means that the flat connections $\widetilde{\Gamma}$ and $\Gamma$ associated
with corresponding symplectic foliations over ${\mathbf{s}}(B)$, are gauge equivalent.
\end{rem}

Now suppose we are given a triple $(E\overset{\pi}{\rightarrow}B,\mathcal{V}%
,\mathbf{s})$ consisting of a fiber bundle over a connected base, a vertical
Poisson tensor $\mathcal{V}$ and a cross-section ${\mathbf{s}}:\,B\rightarrow
E$. Assume that $\mathrm{\operatorname{rank}{\mathcal{V}}=0}$\textrm{\ at
{\textbf{s}}(B)} and
\begin{equation}
{\mathcal{H}}_{V}^{1}(\mathcal{E};{\mathcal{V}})=0
\end{equation}
for a certain open neghborhood $\mathcal{E}$ of ${\mathbf{s}}(B)$ in $E$.
Assume also that there exists a $C^{\infty}(B)$-linear map $\mathbf{m}%
:{\mathcal{X}}_{V}^{\text{Ham}}(\mathcal{E})\rightarrow C^{\infty}%
(\mathcal{E})$ such that
\begin{equation}
Z=\mathcal{V}^{\#}(d\mathbf{m(}Z))\text{.}%
\end{equation}
for every $Z\in{\mathcal{X}}_{V}^{\text{Ham}}(\mathcal{E})$.

We say that a connection $\Gamma$ on $E$ is ${\mathbf{s}}$%
-$\mathit{compatible}$ if condition (3.4) holds. Denote by $C_{B}^{\infty
}(\mathcal{E})$ the subspace of smooth funcions on $\mathcal{E}$ vanishing at
${\mathbf{s}}(B)$. Let $\operatorname{Casim}_{{\mathcal{V}}}^{0}%
(\mathcal{E})\thickapprox\operatorname{Casim}_{{\mathcal{V}}}(\mathcal{E}%
)/\pi^{\ast}C^{\infty}(B)$ be the subspace of Casimir funcions of
$(\mathcal{E}$ ,${\mathcal{V)}}$ vanishing at ${\mathbf{s}}(B)$.

From (3.21), (3.22) we deduce: if $\Gamma$ and $\widetilde{\Gamma}$ are two
${\mathbf{s}}$-compatible Poisson connections on $\mathcal{E}$ ( condition
(2.17) holds), then there exists $\phi_{0}\in\Omega^{1}(B)\otimes
C_{B}^{\infty}(\mathcal{E})$ such that
\begin{equation}
\Gamma-\widetilde{\Gamma}=(\mathcal{V}^{\#}d\phi_{0})^{h}\text{.}%
\end{equation}
Note that $\phi_{0}$ in (3.23) is uniquely determined up to elements from the
space $\Omega^{1}(B)\otimes\operatorname{Casim}_{{\mathcal{V}}}^{0}%
(\mathcal{E})$.

Consider the $C^{\infty}(B)$-module $\mathcal{M}^{k}(\mathcal{E)}=\Omega
^{k}(B)\otimes\operatorname{Casim}_{{\mathcal{V}}}^{0}(\mathcal{E})$. Notice
that the covariant derivative $\partial_{\Gamma}$ associated to a
${\mathbf{s}}$-compatible Poisson connection $\Gamma$ on $\mathcal{E}$ sends
the subspace $\mathcal{M}^{k}(\mathcal{E)}\subset\Omega^{k}(B)\otimes$
$C^{\infty}(\mathcal{E})$ to subspace $\mathcal{M}^{k+1}(\mathcal{E)}%
\subset\Omega^{k+1}(B)\otimes$ $C^{\infty}(\mathcal{E})$. It follows from
(3.21) that there exists a unique operator $\mathcal{\partial}_{0}%
:\mathcal{M}^{k}(\mathcal{E)}\rightarrow\mathcal{M}^{k+1}(\mathcal{E)}$ with
the property: for every ${\mathbf{s}}$-compatible Poisson connection $\Gamma$
on $\mathcal{E}$ the restriction of $\partial_{\Gamma}$ to $\mathcal{M}%
^{k}(\mathcal{E)}$ coincides with $\mathcal{\partial}_{0}$, $\ \partial
_{\Gamma}\mid_{\mathcal{M}^{k}(\mathcal{E)}}=$ $\mathcal{\partial}_{0}$.
Moreover, $\mathcal{\partial}_{0}$ is coboundary operator, $\mathcal{\partial
}_{0}\circ\mathcal{\partial}_{0}=0$.

We say that the germ of $\mathcal{V}$ at ${\mathbf{s}}(B)$ is \textbf{trivial}
if there exists an open neighborhood $\mathcal{E}$ of ${\mathbf{s}}(B)$ in $E$
such that apart from conditions (3.21),(3.22) the second cohomology space of
$\mathcal{\partial}_{0}$ is trivial ,
\begin{equation}
\frac{\ker(\mathcal{\partial}_{0}:\mathcal{M}^{2}(\mathcal{E)}\rightarrow
\mathcal{M}^{3}(\mathcal{E)})}{\operatorname{im}(\mathcal{\partial}%
_{0}:\mathcal{M}^{1}(\mathcal{E)}\rightarrow\mathcal{M}^{2}(\mathcal{E)}%
)}=0\text{.}%
\end{equation}

From Theorem 3.1 we derive the following Poisson analog of the relative
Darboux theorem due to \cite{We1} .

\begin{thm}
Assume that the germ of $\mathcal{V}$ at ${\mathbf{s}}(B)$ is \textit{trivial}%
. Then every two ${\mathbf{s}}$-compatible Poisson tensors $\Pi$ and
$\widetilde{\Pi}$ on $E$ with the same symplectic structure on ${\mathbf{s}%
}(B)$ (condition (3.7)) and the same vertical part
\begin{equation}
\Pi_{V}=\widetilde{\Pi}_{V}\text{ }=\mathcal{V}\text{\ \ \ \ \ \ on
}\mathcal{E}%
\end{equation}
are isomorphic in the sense of (3.11),(3.12).
\end{thm}

\begin{proof}
Let $\Pi$ and $\widetilde{\Pi}$ be two ${\mathbf{s}}$-compatible Poisson
tensors on $E$ satisfying the above hypotheses. Let $(\Gamma,{\mathcal{V}%
},\mathbb{F})$ and $(\widetilde{\Gamma},{\mathcal{V}},\widetilde{\mathbb{F}})$
be the geometric data associated with $\Pi$ and $\widetilde{\Pi}$,
respectavely. Thus, Poisson connections $\Gamma$ and $\widetilde{\Gamma}$ are
${\mathbf{s}}$-compatible and hence (3.23) holds. Pick a $\phi_{0}$ in (3.23)
and define
\begin{equation}
\mathcal{C}:\mathcal{=}\widetilde{\mathbb{F}}-\mathbb{F+}\partial_{\Gamma}%
\phi_{0}+\frac{1}{2}\{\phi_{0}\wedge\phi_{0}\}_{{\mathcal{V}}}\text{.}%
\end{equation}

It follows from (3.23) and the curvature identity (2.19) for $\mathbb{F}$ and
$\widetilde{\mathbb{F}}$ that $\mathcal{C\in M}^{2}(\mathcal{E)}$. Using
(2.18), we deduce: $\mathcal{C}$ is a 2-cocycle,\ $\mathcal{\partial}%
_{0}\mathcal{C=}0$ whose cohomology class does not depend on the choice of
$\phi_{0}$ in (3.26). If this class vanishes, then $\mathcal{C}%
=\mathcal{\partial}_{0}\mathcal{\beta}$ for a $\beta\in\Omega^{1}%
(B)\otimes\operatorname{Casim}_{{\mathcal{V}}}^{0}(\mathcal{E})$ and
$\widetilde{\mathbb{F}}$ and $\mathbb{F}$ satisfy (3.10) for $\phi=\phi
_{0}-\beta$ and $g=\operatorname{id}$.
\end{proof}

It remains to note : for the equivalence of two individual ${\mathbf{s}}%
$-compatible Poisson tensors $\Pi$ and $\widetilde{\Pi}$ instead of (3.24) we
can assume that the cohomology class of the relative 2-cocycle (3.26) is trivial.

\section{Poisson structures from Lie algebroids}

Our goal is to describe a class of connection-dependent coupling tensors on
the dual of the isotropy of a transitive Lie algebroid over a symplectic base.

To begin, we recall some definitions and facts in the theory of Lie algebroids
( for more detail see \cite{Mz,Ku,Va,IKV,CWe} and references given there).

A \textit{Lie algebroid} over a manifold $B$ is a vector bundle $A\to B$
together with a bundle map $\rho:\, A\to TB$, called the \textit{anchor}, and
a Lie algebra structure $\{\,,\,\}_{A}$ on the space $\Gamma(A)$ of smooth
sections of $A$ such that

\begin{enumerate}
\item  For any $a_{1},a_{2}\in\Gamma(A)$,
\begin{equation}
\rho(\{a_{1},a_{2}\}_{A})=[\rho(a_{1}),\rho(a_{2})].
\end{equation}

\item  For any $a_{1},a_{2}\in\Gamma(A)$ and $f\in C^{\infty}(B)$,
\begin{equation}
\{a_{1},fa_{2}\}_{A}=f\{a_{1},a_{2}\}_{A}+(L_{\rho(a_{1})}f)a_{2}.
\end{equation}
\end{enumerate}

The kernel of the anchor $\rho$ is called \textit{isotropy}.

If $A$ and $\widetilde{A}$ are two Lie algebroids over the same base manifold
$B$, then a morphism of Lie algebroids over $B$ is a vector bundle morphism
$\imath:A\rightarrow\widetilde{A}$ over $B$ such that $\widetilde{\rho}%
\circ\imath=\rho$ and such that $\imath(\{a_{1},a_{2}\}_{A})=\{\imath
(a_{1}),\imath(a_{2})\}_{\widetilde{A}}$ for all $a_{1},a_{2}\in\Gamma(A)$. If
$\imath$ is a vector bundle isomorphism we say that $A$ and $\widetilde{A}$
are isomorphic.

A Lie algebroid is called \textit{transitive} if the anchor is a fiberwise surjection.

Let $(A,\rho,[\,,\,]_{A})$ be a transitive Lie algebroid over a connected base
$B$. Then there is an exact sequence of vector bundles
\begin{equation}
\ker\rho\rightarrow A\overset{\rho}{\rightarrow}TB\text{.}%
\end{equation}
It follows from the Lie algebroid axioms that the restriction of the bracket
$\{\,,\,\}_{A}$ to $\Gamma(\ker\rho)$ defines a fiberwise Lie algebra
structure on $\ker\rho$ which will be denoted by $[\,,\,]$. Further,
$(\ker\rho,[\,,\,])$ is a locally trivial Lie algebra bundle with a typical
fiber $\frak{g}$, that is, the structure group of $\ker\rho$ reduces from
$\operatorname{GL}(\frak{g)}$ to the automorphism group $\operatorname{Aut}%
(\frak{g)}$ of the Lie algebra $\frak{g}$ (see \cite{Mz}).

A \textit{connection} on the transitive Lie algebroid $A$ due to Mackenzie
\cite{Mz}, is defined as a right splitting of the exact sequence of vector
bundles (4.3), that is, a vector bundle morphism $\gamma:\,TB\rightarrow A$
such that $\rho\circ\gamma=\operatorname{id}$. Thus $\gamma$ induces
$A=\gamma(TB)\oplus\ker\rho$. The \textit{curvature} of $\gamma$ is the vector
valued $2$-form ${\mathcal{R}}^{\gamma}\in\Omega^{2}(B)\otimes\Gamma(\ker
\rho)$ defined by
\begin{equation}
{\mathcal{R}}^{\gamma}(u_{1},u_{2}):=\{\gamma(u_{1}),\gamma(u_{2}%
)\}_{A}-\gamma([u_{1},u_{2}])
\end{equation}
for $u_{1},u_{2}\in{\mathcal{X}}(B)$.

Given a connection $\gamma$ on $A$, there is a linear Koszul connection
$\nabla^{\gamma}:\,\Gamma(\ker\rho)\rightarrow\Gamma(T^{\ast}B\otimes\ker
\rho)$ on the vector bundle $\ker\rho$, called an \textit{adjoint connection}
\cite{Mz}, and defined by
\begin{equation}
\nabla_{u}^{\gamma}\eta=\{\gamma(u),\eta\}_{A}\qquad(u\in{\mathcal{X}}%
(B),\eta\in\Gamma(\ker\rho)).
\end{equation}
This connection preserves the fiberwise Lie structure on $\ker\rho$,
\begin{equation}
\nabla^{\gamma}([\eta_{1},\eta_{2}])=[\nabla^{\gamma}\eta_{1},\eta_{2}%
]+[\eta_{1},\nabla^{\gamma}\eta_{2}]
\end{equation}
for $\eta_{1},\eta_{2}\in\Gamma(\ker\rho)$. The curvature form
$\operatorname{Curv}^{\nabla^{\gamma}}:\,TB\oplus TB\rightarrow
\operatorname{End}(\ker\rho)$ is given by
\[
\operatorname{Curv}^{\nabla^{\gamma}}(u_{1},u_{2}):=[\nabla_{u_{1}}%
,\nabla_{u_{2}}]-\nabla_{\lbrack u_{1},u_{2}]}%
\]
and related to the curvature of $\gamma$ by the adjoint representation
\begin{equation}
\operatorname{Curv}^{\nabla^{\gamma}}=\operatorname{ad}\circ{\mathcal{R}%
}^{\gamma}.
\end{equation}
Here we use the notation $\operatorname{ad}\circ\eta=[\eta,\cdot]$ for
$\eta\in\Gamma(\ker\rho)$. Furthermore, one can show that ${\mathcal{R}%
}^{\gamma} $ satisfies the Bianchi identity :
\begin{equation}
\underset{(u_{0},u_{1},u_{2})}{\frak{S}}\left(  \nabla_{u_{0}}^{\gamma
}{\mathcal{R}}^{\gamma}(u_{1},u_{2})+{\mathcal{R}}^{\gamma}(u_{0},[u_{1}%
,u_{2}])\right)  =0
\end{equation}
for any $u_{,}u_{1},u_{2}\in{\mathcal{X}}(B)$. Identity (4.8) means that the
$\nabla^{\gamma}$-covariant derivative of ${\mathcal{R}}^{\gamma}$ vanishes,
$\partial^{\nabla^{\gamma}}{\mathcal{R}}^{\gamma}=0$.

\begin{exam}
An important class of transitive Lie algebriods comes from principle bundles.
If we have a $G$-principle bundle $P$ $\overset{\tau}{\longrightarrow}$ $B$,
then there is an exact sequence of of vector bundles
\[
\operatorname{ad}(P)=P\times_{G}\frak{g}\rightarrow TP/G\rightarrow TB
\]
called the Atiyah sequence. Here $\frak{g}$ is the Lie algebra of $G$, $TP/G$
is the quotient manifold with respect to the (right) lifted action to the
cotangent bundle and $\operatorname{ad}(P)$ is the bundle over $B$ associated
with $P $ via the adjoint action of $G$ on $\frak{g}$. The natural isorphism
between smooth sections of $TP/G$ and the space of right invariant vector
fields on $P$ induces the Lie bracket on $\Gamma(TP/G)$. Thus, $A=TP/G$
becomes a transitive Lie algebriod over $B$ whose isotropy is \ the adjoint
bundle $\operatorname{ad}(P)$ (see \cite{Mz,Ku} ). A given principle
connection $\vartheta:TB\rightarrow TP$ with $G$-invariant horizontal
subbundle $\vartheta(TB)$ induces the connection $\gamma:TB\rightarrow TP/G$
on the Lie algebroid $A$. The $\frak{g}$-valued curvature form $\mathcal{K}%
^{\vartheta}\in\Omega^{2}(B;\frak{g)}$ of $\vartheta$ is related with the
curvature ${\mathcal{R}}^{\gamma}:TB\oplus TB\rightarrow P\times_{G}\frak{g}$
by the formula ${\mathcal{R}}^{\gamma}=\operatorname{pr}\circ\tau^{\ast
}\mathcal{K}^{\vartheta}$. Here $\tau^{\ast}\mathcal{K}^{\vartheta}:TP\oplus
TP\rightarrow P\times\frak{g}$ \ is the pull back via the projection $\tau$
and $\operatorname{pr}:P\times\frak{g}\rightarrow P\times_{G}\frak{g}$ is the
natural projection. As is known\ \cite{AM} there are \ ''nonintegrable'' Lie
algebriods, which are transitive \ and can not be realized \ as the Lie
algebroids of principle bundles (also see \cite{Mz,Ku} ).
\end{exam}

\subsection{Connection-dependent coupling tensors}

Let $\nu:\,{\mathcal{N}}\rightarrow B$ be a vector bundle over a connected
symplectic base $(B,\omega)$. Suppose we are given

\begin{itemize}
\item  a transitive Lie algebroid $(A,\rho,\{\,,\,\}_{A})$ over $B$ such that
the isotropy of $A$ coincides with the dual of ${\mathcal{N}}$
\begin{equation}
{\mathcal{N}}^{\ast}\rightarrow A\rightarrow TB,\qquad{\mathcal{N}}^{\ast
}=\ker\rho;
\end{equation}

\item  a connection $\gamma:\,TB\to A$.
\end{itemize}

Recall that ${\mathcal{N}}^{\ast}$ is a Lie algebra bundle with \ fiberwise
Lie algebra structure $[$ , $]$ and typical fiber $\frak{g}$. Hence
${\mathcal{N}}$ can be viewed as a \textit{bundle} of \textit{Lie--Poisson
manifolds} with typical fiber $\frak{g}^{\ast}$.

Denote by $C^{\infty}_{\text{lin}}({\mathcal{N}})$ the space of
\textit{fiberwise linear functions} on ${\mathcal{N}}$. Then we have the
natural identification
\begin{equation}
\ell:\,\Gamma({\mathcal{N}}^{*})\to C^{\infty}_{\text{lin}}({\mathcal{N}})
\end{equation}
given by $\ell(\eta)(x)=\langle\eta(\nu(x)),x\rangle$ for $x\in{\mathcal{N}}$
and $\eta\in\Gamma({\mathcal{N}}^{*})$.

We say that an \textit{Ehresmann connection} on the vector bundle
${\mathcal{N}}$ is \textit{homogeneous} if the horizontal lift of every base
vector field (as a differential operator) preserves the space $C_{\text{lin}%
}^{\infty}({\mathcal{N}})$. Equivalently, the horizontal subbundle is
invariant with respect to dilations $\lambda_{t}:$ ${\mathcal{N}}%
\rightarrow{\mathcal{N}}$ \ ($\lambda_{t}(x)=t\cdot x$, $\ x\in{\mathcal{N}}$,
$t\in\mathbb{R}$). Notice that there is a bijective correspondence between
homogeneous\textit{\ }Ehresmann connections on ${\mathcal{N}}$ and linear
connections\ (covariant derivatives) in the sense of Koszul \cite{GHV}.

Now let us assign to the pair $(A,\gamma)$ a triple $(\Gamma^{A,\gamma
},\Lambda,\mathbb{F}^{A,\gamma})$ consisting of

\begin{itemize}
\item  the homogenious Ehresmann connection $\Gamma^{A,\gamma}$ on
${\mathcal{N}}$ whose horizontal lift is defined by
\begin{equation}
L_{\operatorname{hor}(u)}\varphi=\ell(\{\gamma(u),\ell^{-1}(\varphi)\}_{A})
\end{equation}
for $u\in{\mathcal{X}}(B)$, $\varphi\in C_{\text{lin}}^{\infty}({\mathcal{N}%
})$;

\item  the fiberwise linear vertical Poisson tensor $\Lambda\in\chi_{V}%
^{2}({\mathcal{N}})$ given by
\begin{equation}
\Lambda(d\varphi_{1},d\varphi_{2})=\ell([\ell^{-1}(\varphi_{1}),\ell
^{-1}(\varphi_{2})])
\end{equation}
for $\varphi_{1},\varphi_{2}\in C_{\text{lin}}^{\infty}({\mathcal{N}})$;

\item  the base $2$-form $\mathbb{F}^{A,\gamma}\in\Omega^{2}(B)\otimes
C_{\text{aff}}^{\infty}({\mathcal{N}})$:
\begin{equation}
\mathbb{F}^{A,\gamma}=\omega\otimes1-\ell\circ{\mathcal{R}}^{\gamma}.
\end{equation}
\end{itemize}

For the second term in (4.13) we have $\ell\circ{\mathcal{R}}^{\gamma}%
(u_{1},u_{2})(e)=$ $\left\langle {\mathcal{R}}^{\gamma}(u_{1}(b),u_{2}%
(b)),e\right\rangle $ for $u_{1},u_{2}\in{\mathcal{X}}(B)$ and $e\in
{\mathcal{N}}$, here $b=\nu(e)$. Thus, the homogenious Ehresmann connection
$\Gamma^{A,\gamma}$ is generated by the linear connection on ${\mathcal{N}}$
which is conjugate to the adjoint connection $\nabla^{\gamma} $ in (4.5). The
bivector field $\Lambda$ defines the fiberwise Lie--Poisson structure on the
bundle $\frak{g}^{\ast}\rightarrow{\mathcal{N}}\overset{\nu}{\rightarrow}B$.
The $2$-form $\mathbb{F}^{A,\gamma}$ takes values in the space of fiberwise
affine functions $C_{\text{aff}}^{\infty}({\mathcal{N}})\approx C^{\infty
}(B)\oplus C_{\text{lin}}^{\infty}({\mathcal{N}})$ and includes the base
symplectic $2$-form $\omega$ and the curvature form ${\mathcal{R}}^{\gamma
}:TB\oplus TB\rightarrow{\mathcal{N}}^{\ast}$ in (4.4).

Now we observe that properties (4.6), (4.8) and (4.7) imply relations
(2.17)--(2.19) for $(\Gamma^{A,\gamma},\Lambda,\mathbb{F}^{A,\gamma})$.
Moreover, since $\ell\circ{\mathcal{R}}^{\gamma}\in\Omega^{2}(B)\otimes
C^{\infty}_{\text{lin}}({\mathcal{N}})$, there is a neighborhood $E$ of the
zero section $B\hookrightarrow{\mathcal{N}}$, where the $2$-form
$\mathbb{F}^{A,\gamma}$ is nondegenerate. So applying Theorem~2.1, we arrive
at the following assertion.

\begin{thm}
In a neighborhood $E$ of the zero section $B\hookrightarrow{\mathcal{N}}$ the
transitive Lie algebroid $A$ with a connection $\gamma$ induces a coupling
tensor $\Pi^{A,\gamma}$ associated with the geometric data $(\Gamma^{A,\gamma
},\Lambda,\mathbb{F}^{A,\gamma})$ in $\emph{(4.11)-(4.13)}$. If the kernel
$\ker{\mathcal{R}}^{\gamma}\subset TB$ of the curvature $2$-form
${\mathcal{R}}^{\gamma}$ is a coisotropic distribution with respect to the
base symplectic form $\omega$, then the coupling tensor $\Pi^{A,\gamma}$ is
well-defined on the entire total space ${\mathcal{N}}$.
\end{thm}

To justify the second part of Theorem~4.1, let us consider the coordinate
representation for $\Pi^{A,\gamma}$.

Let $(\xi,x)=(\xi^{1},\dots,\xi^{2k};x^{1},\dots,x^{r})$ be a (local)
coordinate system on ${\mathcal{N}}$, where $(\xi^{i})$ are coordinates on the
base $B$ and $(x^{\sigma})$ are coordinates on the fibers of ${\mathcal{N}}$
associated with a basis of local sections $(X_{\sigma})$. Then we have

\begin{itemize}
\item  the \textit{symplectic form on the base\/}: $\omega=\frac{1}{2}%
\sum_{i,j}\omega_{ij}(\xi)\,d\xi^{i}\wedge d\xi^{j},$ $\ \omega^{is}%
\omega_{sj}=\delta_{j}^{i};$

\item  the \textit{curvature form\/}: ${\mathcal{R}}^{\gamma}=\frac{1}{2}%
\sum_{i,j,\sigma}{\mathcal{R}}_{ij\sigma}(\xi)d\xi^{i}\wedge d\xi^{j}\otimes
dx^{\sigma};$

\item  the \textit{connection form\/}: $\Gamma^{A,\gamma}=\sum_{i,\sigma
}\Gamma_{i}^{\sigma}d\xi^{i}\otimes\frac{\partial}{\partial x^{\sigma}}%
,\qquad\Gamma_{i}^{\sigma}=\Gamma_{i\sigma^{\prime}}^{\sigma}(\xi
)x^{\sigma^{\prime}};$

\item  the \textit{base 2-form\/} (4.13): $\mathbb{F}^{A,\gamma}=\frac{1}%
{2}\sum_{i,j}d\xi^{i}\wedge d\xi^{j}\otimes F_{ij},$ where
\begin{equation}
F_{ij}=\omega_{ij}-\sum_{\sigma}{\mathcal{R}}_{ij\sigma}x^{\sigma}.
\end{equation}
\end{itemize}

Let $(\eta^{\sigma})$ be the dual basis of local sectons of ${\mathcal{N}%
}^{\ast}$, $\left\langle \eta^{\sigma},X_{\sigma^{^{\prime}}}\right\rangle
=\delta_{\sigma^{^{\prime}}}^{\sigma}$. Then \ with respect to the induced
basis of local sections $(\Xi_{i}=\gamma(\frac{\partial}{\partial\xi^{i}}),$
$\eta^{\sigma})$ of \ $A$ the Lie algebroid structure takes the form :
\[
\{\Xi_{i},\Xi_{j}\}_{A}=\sum_{\nu}{\mathcal{R}}_{ij\sigma}\eta^{\nu},\text{
\ \ }\{\Xi_{i},\eta^{\sigma}\}_{A}=-\sum_{\nu}\Gamma_{i\nu}^{\sigma}\eta^{\nu
},\text{ \ \ \ }\{\eta^{\sigma},\eta^{\sigma^{^{\prime}}}\}_{A}=\sum_{\nu
}\lambda_{\nu}^{\sigma\sigma}\eta^{\nu}\text{.}%
\]
Consider the open domain containing the zero section $B=\{x^{1}=0,\dots$,
$x^{r}=0\}$:
\begin{equation}
E=\{(\xi,x)\in{\mathcal{N}}\mid\det(\!(\omega_{ij}-\sum_{\sigma}{\mathcal{R}%
}_{ij\sigma}x^{\sigma})\!)\neq0\}.
\end{equation}
Then the coupling tensor $\Pi^{A,\gamma}$ is well-defined on $E$ and has the
representation
\begin{equation}
\Pi^{A,\gamma}=\frac{1}{2}\sum_{i,j}H^{ij}(\xi,x)\operatorname{hor}%
(\partial_{i})\wedge\operatorname{hor}(\partial_{j})+\frac{1}{2}\sum
_{\sigma\sigma^{\prime}}\Lambda^{\sigma\sigma^{\prime}}(\xi,x)\frac{\partial
}{\partial x^{\sigma}}\wedge\frac{\partial}{\partial x^{\sigma^{\prime}}}.
\end{equation}
Here $\operatorname{hor}(\partial_{i})=\partial/\partial\xi^{i}-\sum_{\sigma
}\Gamma_{i}^{\sigma}\partial/\partial x^{\sigma}$ and the matrix functions
$(\!(H^{ij})\!)$ and $(\!(\Lambda^{\sigma\sigma^{\prime}})\!)$ are defined by
\begin{equation}
\sum_{s}H^{is}F_{sj}=-\delta_{j}^{i},\qquad\Lambda^{\sigma\sigma^{\prime}%
}=\sum_{\nu}\lambda_{\nu}^{\sigma\sigma^{\prime}}(\xi)x^{\nu},
\end{equation}
where $\lambda_{\nu}^{\sigma\sigma^{\prime}}(\xi)$ are the structure constants
of the Lie algebra ${\mathcal{N}}_{\xi}\approx\frak{g}$.

The Poisson brackets of the coupling tensor $\Pi^{A,\gamma}$ on the domain
(4.15) take the form:
\begin{align}
\{\xi^{i},\xi^{j}\} &  =H^{ij}=(-\omega^{ij}+\omega^{ii^{\prime}}{\mathcal{R}%
}_{i^{\prime}j^{\prime}\sigma}\omega^{j^{\prime}j}x^{\sigma})+O_{2}%
,\nonumber\\
\{\xi^{i},x^{\sigma}\} &  =-H^{is}\Gamma_{s}^{\sigma}=\omega^{is}%
\Gamma_{s\sigma^{\prime}}^{\sigma}x^{\sigma^{\prime}}+O_{2},\\
\{x^{\sigma},x^{\sigma^{\prime}}\} &  =\Lambda^{\sigma\sigma^{\prime}}%
+H^{ij}\Gamma_{i}^{\sigma}\Gamma_{j}^{\sigma^{\prime}}=(\lambda_{\nu}%
^{\sigma\sigma^{\prime}}x^{\nu}-\omega^{ij}\Gamma_{i\nu}^{\sigma}\Gamma
_{j\nu^{\prime}}^{\sigma^{\prime}}x^{\nu}x^{\nu^{\prime}})+O_{3}.\nonumber
\end{align}
Here the summation is taken with respect to repeated indices and $O_{k}$
denotes a term having zero of order~$k$ at every point in $B$.

Finally, using standard facts from linear symplectic geometry, it is easy to
show that under the coisotropic hypothesis for $\ker{\mathcal{R}}^{\gamma}$,
the matrix $(\!(F^{ij})\!)$ in (4.14) is totally nondegenerate and hence
domain (4.15) coincides with the total space ${\mathcal{N}}$.

\begin{rem}
In the case when $A$ is the coadjoint bundle of a principle bundle, connection
dependent Poisson structures of type $\Pi^{A,\gamma}$ were studied in
\cite{MoMR,Mo}.
\end{rem}

\begin{exam}
\textrm{Suppose we are given a vector bundle $\nu:\,{\mathcal{L}}\to Q$
equipped with }

\begin{itemize}
\item \textrm{a fiberwise Lie algebra structure \ }$[\eta^{\sigma}%
,\eta^{\sigma^{\prime}}]_{{\mathcal{L}}}=\sum_{\nu}\lambda_{\nu}^{\sigma
\sigma^{\prime}}(q)\eta^{\nu},$\textrm{\ }

\item \textrm{a linear connection \ }$\nabla_{\partial/\partial q_{i}}%
\eta^{\sigma}=-\sum_{\sigma^{\prime}}\theta_{i\sigma^{\prime}}^{\sigma}%
(q)\eta^{\sigma^{\prime}}.$\textrm{\ }
\end{itemize}

\textrm{Here $(\eta^{\sigma})$ is a basis of local sections of ${\mathcal{L}}$
and $q=(q^{i})$ are local coordinates on the base $Q$. Assume that }

\textrm{(i) $\nabla$ preserves $[\,,\,]_{{\mathcal{L}}}$ (condition (4.6)); }

\textrm{(ii) there exists a vector bundle morphism ${\mathcal{R}}:TQ\times
TQ\rightarrow{\mathcal{L}}$
\[
{\mathcal{R}}\bigg(\frac{\partial}{\partial q^{i}},\frac{\partial}{\partial
q^{j}}\bigg)=\sum_{\nu}{\mathcal{R}}_{ij\nu}(q)\eta^{\nu}%
\]
which is related to the curvature $2$-form $\operatorname{Curv}^{\nabla}$ on
$Q$ by formula (4.7); }

\textrm{(iii) ${\mathcal{R}}$ satisfies the modified Bianchi identity (4.8). }

Then the triple $(\nabla,$\textrm{${\mathcal{R}}$,\ }$[,]_{{\mathcal{L}}})$
defines the transitive Lie algebroid on $A=TQ\oplus$\textrm{${\mathcal{L}}$ (
\cite{Mz} ) such that }$\operatorname{pr}_{1}:TQ\oplus$%
\textrm{${\mathcal{L\rightarrow}}$}$TQ$ is the anchor, \textrm{\ ${\mathcal{L}%
}$ is the isotropy , }$\nabla$ and \textrm{${\mathcal{R}}$ is the adjoint
connection and the curvature of the connection }$\mathcal{\gamma}%
_{0}:TQ\rightarrow TQ\oplus$\textrm{${\mathcal{L}}$ (canonical injection).
Consider the pull back\ \ }$\widetilde{A}\rightarrow T^{\ast}Q$ of $A$ via the
natural projection $T^{\ast}Q\rightarrow Q$. Denote also by $(\widetilde
{\nabla}$, $\widetilde{\mathrm{{\mathcal{R}}}}$\textrm{\ ,}$[,]_{\widetilde
{\mathrm{{\mathcal{L}}}}})$ the cotangent pull back of the original triple
$(\nabla,$\textrm{${\mathcal{R}}$,\ }$[,]_{{\mathcal{L}}})$ and by
$\widetilde{\mathrm{{\mathcal{L}}}}\mathrm{\ }\rightarrow TQ$ the pull back of
the bundle \textrm{${\mathcal{L\rightarrow}}Q$. }Consider the
\textrm{canonical symplectic structure }$\omega=\sum_{i}dp^{i}\wedge
dq^{i}=dp\wedge dq$ on \textrm{$T^{\ast}Q$ . Then the triple }$(\widetilde
{\nabla}$, $\widetilde{\mathrm{{\mathcal{R}}}}$, \textrm{\ }$[,]_{\widetilde
{\mathrm{{\mathcal{L}}}}})$ induces a transitive Lie algebroid on
$\widetilde{A}$ over the symplectic base \textrm{$(B=T^{\ast}Q,\omega=dp\wedge
dq)$ (this is an inverse-image algebroid\textbf{\ \cite{Mz,Ku}). }Moreover,
}$\widetilde{\mathrm{{\mathcal{L}}}}$ is the isotropy of $\widetilde{A}$ and
the pull back $\widetilde{\mathcal{\gamma}_{0}}$ $:T($\textrm{$T^{\ast
}Q)\rightarrow$}$\widetilde{A}$ is the connection on $\widetilde{A}$ whose
curvature is just $\widetilde{\mathrm{{\mathcal{R}}}}$.\textrm{\ Thus, the
kernel of } $\widetilde{\mathrm{{\mathcal{R}}}}$\ \textrm{is a Lagrangian
distribution on }$T^{\ast}Q$\textrm{\ with respect to the form $dp\wedge dq$.
Hence the coupling tensor associated with the pair $(\widetilde{A}%
,\widetilde{\mathcal{\gamma}_{0}})$ is well defined on the entire total space
of the dual }$\widetilde{\mathrm{{\mathcal{L}}}}^{\ast}$\textrm{and the
Poisson bracket in (4.18) takes the following coordinate form}

\textrm{\ } \textrm{\
\begin{align*}
\{p^{i},p^{j}\} &  ={\mathcal{R}}_{ij\nu}(q)x^{\nu},\text{ \ \ \ }%
\{p^{i},q^{j}\}=\delta^{ij},\qquad\{q^{i},q^{j}\}=0,\\
\{p^{i},x^{\sigma}\} &  =-\theta_{i\sigma^{\prime}}^{\sigma}(q)x^{\sigma
^{\prime}},\text{ \ \ }\{p^{i},x^{\sigma}\}=\{q^{i},x^{\sigma}\}=0,\\
\{x^{\sigma},x^{\sigma^{\prime}}\} &  =\lambda_{\nu}^{\sigma\sigma^{\prime}%
}(q)x^{\nu}.
\end{align*}
}

On the other hand, it is of interest to note: \textrm{this Poisson structure
coincides with the Courant structure \cite{Co} on the dual }$A^{\ast}%
=TQ^{\ast}\oplus$\textrm{${\mathcal{L}}$}$^{\ast}$ of the Lie algebroid $A$.
\textrm{Notice also that such a type of Poisson structures arises from the
study of Hamiltonian structures for Wong's equations \cite{MoMR,Mo,La}. }
\end{exam}

\subsection{Varying the connection and the Lie algebroid structure}

Let us address the following question: how does the coupling tensor
$\Pi^{A,\gamma}$ defined in Theorem~4.1 depend on the choice of the connection
$\gamma$ and the Lie algebroid structure on $A$? We will investigate this
issue in two steps. Let $A$ be a transitive Lie algebroid over a connected
symplectic base $(B,\omega)$. Let ${\mathcal{L}}$ be the isotropy of $A$ and
let ${\mathcal{N}}={\mathcal{L}}^{*}$ be the dual. The fiberwise Lie structure
on ${\mathcal{L}}$ will be denoted by $[\,,\,]_{{\mathcal{L}}}$.

\textbf{I}. Suppose that we have two connections on $A$:
\[
\gamma:\,TB\rightarrow A\quad\text{and}\quad\widetilde{\gamma}:\,TB\rightarrow
A.
\]
Consider adjoint connections and curvature forms $\nabla^{\gamma}%
,{\mathcal{R}}^{\gamma}$ and $\nabla^{\tilde{\gamma}},{\mathcal{R}}%
^{\tilde{\gamma}}$ associated to $\gamma$ and $\tilde{\gamma}$ respectively.
There is a vector bundle map $\mu:\,TB\rightarrow{\mathcal{L}}$ such that
\[
\tilde{\gamma}(u)=\gamma(u)+\mu(u)\qquad\text{for}\quad u\in{\mathcal{X}}(B).
\]
We can think of $\mu$ as a ${\mathcal{L}}$-valued $1$-form on $B$, $\mu
\in\Omega^{1}(B)\otimes\Gamma({\mathcal{L}})$. Then we have \cite{Mz}:
\begin{align}
\nabla_{u}^{\tilde{\gamma}} &  =\nabla_{u}^{\gamma}+\operatorname{ad}\circ
\mu(u),\qquad u\in{\mathcal{X}}(B),\\
{\mathcal{R}}^{\tilde{\gamma}} &  ={\mathcal{R}}^{\gamma}+\partial
_{\nabla^{\gamma}}\mu+\frac{1}{2}[\mu\wedge\mu]_{{\mathcal{L}}}.
\end{align}
Here $\partial_{\nabla^{\gamma}}:\,\Omega^{k}(B)\otimes\Gamma({\mathcal{L}%
})\rightarrow\Omega^{k+1}(B)\otimes\Gamma({\mathcal{L}})$ is the covariant
exterior derivative associated with the linear connection $\nabla^{\gamma}$,
and in the last term in (4.20) we use the standard bracket on the graded
algebra of $L$-valued forms on $B$ generated by the fiberwise Lie algebra
structure $[\,,\,]_{{\mathcal{L}}}$. Now let us consider the geometric data
$(\Gamma^{A,\gamma},\Lambda,\mathbb{F}^{A,\gamma})$ and $(\Gamma
^{A,\tilde{\gamma}},\Lambda,\mathbb{F}^{A,\tilde{\gamma}})$ defined in (4.11)--(4.13).

It follows from (4.19), (4.20) that $\Gamma^{A,\tilde\gamma},\Gamma^{A,\gamma
}$ and $\mathbb{F}^{A,\tilde\gamma},\mathbb{F}^{A,\gamma}$ satisfy relations
(3.9), (3.10) for $g=\operatorname{id}$, ${\mathcal{V}}=\Lambda$, and
\begin{equation}
\phi=\ell\circ\mu\in\Omega^{1}(B)\otimes C^{\infty}_{\text{lin}}({\mathcal{N}%
}).
\end{equation}
Thus the geometric data $(\Gamma^{A,\gamma},\Lambda,\mathbb{F}^{A,\gamma})$
and $(\Gamma^{A,\tilde\gamma},\Lambda,\mathbb{F}^{A,\tilde\gamma})$ are
equivalent. Consider the corresponding coupling tensors $\Pi^{A,\gamma}$,
$\Pi^{A,\tilde\gamma}$ on ${\mathcal{N}}$. Then the zero section
$B\hookrightarrow{\mathcal{N}}$ with a given symplectic form $\omega$ is a
common symplectic leaf of $\Pi^{A,\gamma}$ and $\Pi^{A,\tilde\gamma}$. So, we
can apply to $\Pi^{A,\gamma}$ and $\Pi^{A,\tilde\gamma}$ the neighborhood
equivalence Theorem~3.1.

\begin{prop}
Coupling tensors $\Pi^{A,\gamma}$ and $\Pi^{A,\tilde{\gamma}}$ associated with
arbitrary connections $\gamma$ and $\tilde{\gamma}$ on $A$ are isomorphic over
$B$, that is, there are open neighborhoods ${\mathcal{O}}$, $\widetilde
{{\mathcal{O}}}$ of the zero section $B\hookrightarrow{\mathcal{N}}$ and a
diffeomorphism ${\mathbf{f}}:\,{\mathcal{O}}\rightarrow\widetilde
{{\mathcal{O}}}$ identical on $B$ such that ${\mathbf{f}}^{\ast}\Pi
^{A,\tilde{\gamma}}=\Pi^{A,\gamma}$.
\end{prop}

The equivalence class of isomorphic Poisson structures $\Pi^{A,\gamma}$ will
be called an $\omega$-\textbf{coupling structure} of a transitive Lie
algebroid $A$.

\textbf{II}. Let $A$ and $\widetilde{A}$ be two transitive Lie algebroids over
the same connected base $(B,\omega)$. Assume that $A$ and $\widetilde{A}$ are
isomorphic and $\imath:\widetilde{A}\rightarrow\,A$ is a Lie algebroid
isomorphism. Without loss of generality, we can also assume that
\begin{align}
A &  =TB\oplus{\mathcal{L}},\\
\widetilde{A} &  =TB\oplus\widetilde{{\mathcal{L}}},
\end{align}
and the corresponding anchors $\rho:\,A\rightarrow TB$, $\widetilde{\rho
}:\,\widetilde{A}\rightarrow TB$ coincide with the canonical projections
$\rho=\operatorname{pr}_{1}$, $\widetilde{\rho}=\widetilde{\operatorname{pr}%
}_{1}$. It is clear that the restriction
\begin{equation}
g=\imath\big|_{\widetilde{{\mathcal{L}}}}:\,\widetilde{{\mathcal{L}}%
}\rightarrow{\mathcal{L}}%
\end{equation}
is a vector bundle isomorphism preserving the fiberwise Lie algebra structure
on ${\mathcal{L}}$ and $\widetilde{{\mathcal{L}}}$. We observe that $\imath$
takes an element $u\oplus\eta$ in $\,\widetilde{A}$ into the element
$\imath(u\oplus\eta)$ in $A$ of the form
\begin{equation}
\imath(u\oplus\eta)=u\oplus(g(\eta)+\mu(u)),
\end{equation}
where $\mu:\,TB\rightarrow{\mathcal{L}}$ is a vector bundle morphism. Thus,
$\imath$ is characterized by the pair $(g,\mu)$. Define connections
$\gamma_{0}$ on $A$ and $\tilde{\gamma}_{0}$ on $\widetilde{A}$ as the
canonical injections:
\begin{align}
u &  \mapsto\gamma_{0}(u)=u\oplus0\in TB\oplus{\mathcal{L}},\\
u &  \mapsto\tilde{\gamma}_{0}(u)=u\oplus0\in TB\oplus\widetilde{{\mathcal{L}%
}}.
\end{align}
Then we get
\begin{align}
g([a_{1},a_{2}]_{{\mathcal{L}}}) &  =[g(a_{1}),g(a_{2})]_{\widetilde
{{\mathcal{L}}}}\qquad(a_{1},a_{2}\in\Gamma({\mathcal{L}})),\\
g\circ\nabla_{u}^{\tilde{\gamma}_{0}}\circ g^{-1} &  =\nabla_{u}%
+\operatorname{ad}\circ\mu(u)\qquad(u\in{\mathcal{X}}(B)),\\
g\circ{\mathcal{R}}^{\tilde{\gamma}_{0}} &  ={\mathcal{R}}^{\gamma_{0}%
}+\partial_{\nabla_{\gamma_{0}}}\mu+\frac{1}{2}[\mu\wedge\mu]_{{\mathcal{L}}}.
\end{align}
Relations (4.28)--(4.30) lead to the equivalence relations (3.8)--(3.10) for
geometric data $(\Gamma^{\gamma_{0}},\Lambda,\mathbb{F}^{A,\gamma_{0}})$ and
$(\Gamma^{\tilde{\gamma}_{0}},\widetilde{\Lambda},\mathbb{F}^{\widetilde
{A},\tilde{\gamma}_{0}})$ associated to pairs $(A,\gamma_{0})$ and
$(\widetilde{A},\tilde{\gamma}_{0})$, respectively. As a consequence of
Theorem~3.1, we get the proposition.

\begin{prop}
There is the neighborhood equivalence between coupling tensors $\Pi
^{A,\gamma_{0}}$ and $\Pi^{\widetilde{A},\tilde{\gamma}_{0}}$.
\end{prop}

Finally, combining Proposition~4.1 with Proposition~4.2, we obtain the main result.

\begin{thm}
Let $A$ and $\widetilde{A}$ be two transitive Lie algebroids over the same
connected symplectic base $(B,\omega)$, and let $\gamma:\,TB\rightarrow A$,
$\tilde{\gamma}:\,TB\rightarrow\widetilde{A}$ be two connections. Consider
coupling tensors $\Pi^{A,\gamma}$ and $\Pi^{\widetilde{A},\tilde{\gamma}}$
associated to $(A,\gamma)$ and $(\widetilde{A},\tilde{\gamma})$, respectively.

(i) Assume that $A$ is isomorphic to $\widetilde{A}$. Then under the arbitrary
choice of connections $\gamma$, $\tilde{\gamma}$, there exists a
diffeomorphism ${\mathbf{f}};\,{\mathcal{O}}\rightarrow\widetilde
{{\mathcal{O}}}$ from a neighborhood ${\mathcal{O}}$ of the zero section
$B\hookrightarrow{\mathcal{N}}={\mathcal{L}}^{\ast}$ \emph{(}${\mathcal{L}}$
is the isotropy of $A$\emph{)} onto a neighborhood $\widetilde{{\mathcal{O}}}
$ of the zero section $B\hookrightarrow\widetilde{{\mathcal{N}}}%
=\widetilde{{\mathcal{L}}}^{\ast}$ \emph{(}$\widetilde{{\mathcal{L}}}$ is the
isotropy of $\widetilde{A}$\emph{)} such that ${\mathbf{f}}\big|%
_{B}=\operatorname{id}_{B}$ and
\begin{equation}
{\mathbf{f}}^{\ast}\Pi^{\widetilde{A},\tilde{\gamma}}=\Pi^{A,\gamma\text{
\ \ \ }}\text{and\ \ \ }{\mathbf{f}}\big|_{B}=\operatorname{id}_{B}\text{.}%
\end{equation}

(ii) On the contrary, the equivalence between coupling tensors $\Pi^{A,\gamma
}$ and $\Pi^{\widetilde{A},\tilde{\gamma}}$ (in the sense of (4.31)) implies
the isomorphism between the corresponding Lie algebroids $A$\ and
$\widetilde{\text{ }A}$.
\end{thm}

Now suppose we start with some data $({\mathcal{L}},[$ , $]_{{\mathcal{L}}}%
$,$\frak{g})$, where $({\mathcal{L}},[$ , $]_{{\mathcal{L}}})$ ia locally
trivial bundle of Lie algebras over a coonected symplectic base $(B,\omega)$,
$\frak{g}$ is the typical fiber. Let $\nabla$ be a linear connection in
${\mathcal{L}}$ preserving the fiberwise Lie algebra structure $[$ ,
$]_{{\mathcal{L}}}$ (condition (4.6)) and ${\mathcal{R\in}}\Omega
^{2}(B)\otimes\Gamma({\mathcal{L}})$ be a vector valued 2-form which is
compatible with $(\nabla$,$[$ , $]_{{\mathcal{L}}})$ by means of (4.7) and
(4.8). In this case, we say that the pair $(\nabla,\mathcal{R})$ is admissible
for $[$ , $]_{{\mathcal{L}}}$. Accoding to \cite{Mz} the pair $(\nabla
,\mathcal{R})$ induces a unique transitive Lie algebroid structure $\{$ ,
$\}_{\nabla,\mathcal{R}}$ on $A$ $=TB\oplus{\mathcal{L}}$ such that the anchor
is the natural projection, $({\mathcal{L}},[,]_{{\mathcal{L}}})$ is the
isotropy,$\nabla$ is the adjoint connection associated with connection
$\gamma^{0\text{ }}$ in (4.26) and $\mathcal{R}$ is the curvature of
$\gamma^{0}$. The coupling tensor on ${\mathcal{L}}$ associated to $\{$ ,
$\}_{\nabla,\mathcal{R}}$ and $\gamma^{0}$ will be denoted by $\Pi
^{\nabla,{\mathcal{R}}}$. Consider the subbundle $\operatorname{Cent}%
({\mathcal{L)\subset L}}$ whose typical fiber is the center of the Lie algebra
$\frak{g}$. Then $\operatorname{Cent}({\mathcal{L)}}$ is invariant with
respect to the connection $\nabla$ and the restriction $\nabla_{0}=\nabla
\mid_{\operatorname{Cent}({\mathcal{L)}}}$is a flat connection which does not
depend on the choice of $\nabla$ in the class of adjoint connections of the
Lie algebroid ( see \cite{IKV} ). Thus the covariant derivative $\partial
_{0}:\Omega^{k}(B;\operatorname{Cent}({\mathcal{L))\rightarrow}}\Omega
^{k+1}(B;\operatorname{Cent}({\mathcal{L))}}$ associated with $\nabla_{0}$ is
a coboundary operator , $\partial_{0}\circ\partial_{0}=0$. Notice that the
comology of $\partial_{0}$ coincides with the cohomology of the
\textit{abelian} Lie subalgebroid $A_{0}=TB\oplus\operatorname{Cent}%
({\mathcal{L)}}$ in $(A,\{,\}_{\nabla,\mathcal{R}})$.

Let $(\widetilde{\nabla\text{,}}\widetilde{{\mathcal{R}}})$ be a second
admissible pair for $[$ , $]_{{\mathcal{L}}}$ and $(A,\{,\}_{\widetilde
{\nabla\text{,}}\widetilde{{\mathcal{R}}}})$ be the corresponding Lie
algebroid. Assume that connection $\widetilde{\nabla}$ and $\nabla$ on
${\mathcal{L}}$ are related by (4.19) for a certain $\mu\in\Omega
^{1}(B)\otimes\Gamma({\mathcal{L}})$. This condition means that the structures
of abelian Lie algebroids on $A_{0}$ coming from the brackets $\ \ \{$ ,
$\}_{\nabla,\mathcal{R}}$ and $\{$ , $\}_{\widetilde{\nabla\text{,}}%
\widetilde{{\mathcal{R}}}}$ \ coincide. It follows from (4.8) and (4.19) that
\begin{equation}
\mathcal{C}:=\widetilde{{\mathcal{R}}}-{\mathcal{R}}-\partial_{\nabla}%
\mu-\frac{1}{2}[\mu\wedge\mu]_{{\mathcal{L}}}\text{.}%
\end{equation}
is a 2-cocyle $\mathcal{C\in}\Omega^{2}(B)\otimes\operatorname{Cent}%
({\mathcal{L)}}$, $\partial_{0}\mathcal{C}=0$ whose cohomology class does not
depend on the choice of $\mu$ in (4.19). Moreover we observe: the Lie
algebroid structures $\{$ , $\}_{\nabla,\mathcal{R}}$ and $\{$ ,
$\}_{\widetilde{\nabla\text{,}}\widetilde{{\mathcal{R}}}}$ are isomorphic if
and only if $[\mathcal{C]}=0$. Then as a consenquence of Theorem 4.2., we get
the following ''linear'' analog of Theorem 3.2.

\begin{prop}
Under assuption (4.19) the coupling tensors $\Pi^{\nabla,{\mathcal{R}}}$ and
$\Pi^{\widetilde{\nabla},\widetilde{{\mathcal{R}}}\text{ }}$ are isomorphic
over $B$ if and only if the cohomology class of the relative 2-cocycle
$\mathcal{C}$ in (4.32) is zero. In particular, this is true in the case when
the second cohomology space of the abelian Lie algebroid $A_{0} $ is trivial.
\end{prop}

\begin{rem}
Assume that the typical fiber $\frak{g}$ is \textit{reductive, }that is,
$\frak{g}$ $=\operatorname{Cent}(\frak{g}{\mathcal{)}}\oplus\lbrack
\frak{g},\frak{g]}$, where $[\frak{g},\frak{g]}$ is a semisimple Lie algebra.
Then vanishing of the second cohomology of $A_{0}$ leads to the same property
for the second cohomology of the transitive Lie algebroid $A$, $\mathcal{H}%
^{2}(A)=0$\cite{IKV}. This condition appears also under the study of the
formal Poisson equivalence \cite{IKV}.
\end{rem}

\section{Linearized Poisson models over a single\newline symplectic leaf}

In this section we will show that for every Poisson manifolds with a given
closed symplectic leaf $B$ there is a well defined notion of a
\textit{linearized Poisson structure} at $B$. This linearized structure is
defined as an equivalence class of isomorphic Poisson structures which live
naturally on the normal bundle to the symplectic leaf $B$. In the
zero-dimensional case ($\dim B=0$), our definition coincides with the notion
of a linear approximation of a Poisson structure at a point of
$\operatorname{rank}0$ arising in the context of the linearization problem
\cite{We4}.

\subsection{First approximations}

Let $(M,\Psi)$ be a Poisson manifolds equipped with a Poisson bracket
\begin{equation}
\{F,G\}=\Psi(dF,dG).
\end{equation}
Suppose that we are given a closed (embedded) symplectic leaf $(B,\omega)$ of
$M$ with symplectic structure $\omega$. Consider the \textit{normal bundle} to
the symplectic leaf $B$:
\begin{equation}
{\mathcal{N}}=T_{B}M/TB.
\end{equation}

The well known fact is that the original Poisson structure on $M$ induces a
fiberwise Lie--Poisson structure on the normal bundle ${\mathcal{N}}$ which is
given by the vertical Poisson bivector field $\Lambda\in\chi^{2}({\mathcal{N}%
})$ called a \textit{linearized transverse Poisson structure} of the leaf $B$
\cite{We4}. At each fiber $N_{b}$ over $b\in B$ the Lie--Poisson structure
$\Lambda_{b}\in\chi^{2}({\mathcal{N}}_{b}))$ can be defined as the
linearization of the transverse Poisson structure at $b$ due to the splitting
theorem. To compare the original Poisson tensor $\Psi$ with $\Lambda$, it is
natural to consider a pull back of $\Psi$ onto ${\mathcal{N}}$ via an
exponential map.

By an \textbf{exponential map}, we mean a diffeomorphism ${\mathbf{f}%
}:\,{\mathcal{N}}\to M$ from the normal bundle ${\mathcal{N}}$ onto a tubular
neighborhood of the leaf $B$ in $M$ such that

(i) ${\mathbf{f}}$ is compatible with the zero section $s_{0}%
:\,B\hookrightarrow{\mathcal{N}}$, that is, ${\mathbf{f}}\circ s_{0}=s_{0}$; and

(ii) the composite map
\[
{\mathcal{N}}_{b}\hookrightarrow T_{b}({\mathcal{N}})\xrightarrow{d_b\fff
}T_{b}M\overset{\nu_{b}}{\rightarrow}\mathcal{N}_{b}%
\]
is the identity. Here the last mapping is the canonical projection $\nu
:T_{B}M\rightarrow T_{B}M/TB$.

It follows from the tubular neighborhood theorem that an exponential map
always exists \cite{LMr}. By Proposition 3.1 we deduce the following statement.

\begin{prop}
Let ${\mathbf{f}}^{\ast}\Psi\in\chi^{2}({\mathcal{N}})$ be the pull back of
the Poisson tensor $\Psi$ via an exponential map ${\mathbf{f}}$. Then the zero
section $B\hookrightarrow{\mathcal{N}}$ is a closed symplectic leaf of
${\mathbf{f}}^{\ast}\Psi$ with symplectic structure $\omega$. Moreover, there
exists an open neighborhood $E$ of $B$ in ${\mathcal{N}}$ such that
${\mathbf{f}}^{\ast}\Psi$ is a coupling tensor on $E$. For the vertical tensor
$({\mathbf{f}}^{\ast}\Psi)_{V}$ defined in \emph{(3.1),} we have
\begin{equation}
({\mathbf{f}}^{\ast}\Psi)_{V}=\Lambda+O_{2}\qquad\text{on}\quad E,
\end{equation}
that is, the linearized transverse Poisson structure $\Lambda$ gives a linear
approximation to the vertical part of ${\mathbf{f}}^{\ast}\Psi$.
\end{prop}

\begin{dfn}
\textrm{A $0$-section compatible Poisson tensor $\Pi$ defined (as a coupling
tensor) on an open (tubular) neighborhood $E$ of $B$ in ${\mathcal{N}}$ is
said to be a \textbf{first approximation} to $\Psi$ at the leaf $B$ if }

\textrm{(i) the intrinsic Ehresmann connection $\Gamma$ (2.4) of $\Pi$ is
\textit{homogeneous} on $E$; }

\textrm{(ii) the vertical part $\Pi_{V}$ in (3.1) coincides with the
linearized transverse Poisson structure $\Lambda$ of $B$,
\begin{equation}
\Pi_{V}=\Lambda\qquad\text{on}\quad E;
\end{equation}
}

\textrm{(iii) there exists an exponential map ${\mathbf{f}}:\,{\mathcal{N}%
}\rightarrow M$ such that
\begin{equation}
{\mathbf{f}}^{\ast}\Psi=\Pi+O_{2}\qquad\text{on}\quad E.
\end{equation}
}
\end{dfn}

\begin{thm}
Let $(M,\Psi,B,\omega)$ be a Poisson manifold with a closed symplectic leaf
$(B,\omega)$. Then for a given exponential map ${\mathbf{f}}$ there exists a
unique first approximation $\Pi^{{\mathbf{f}}}$ to $\Psi$ at $B$ satisfying
\emph{(5.5)}. The Poisson bivector field $\Pi^{{\mathbf{f}}}$ does not depend
on the choice of ${\mathbf{f}}$ up to $0$-section neighborhood isomorphism.
\end{thm}

\begin{dfn}
\textrm{The equivalence class of isomorphic Poisson tensors $\Pi^{{\mathbf{f}%
}}$ is said to be the \textbf{linearized Poisson structure} of the leaf $B$. }

\begin{rem}
If the symplectic leaf $B$ is not closed, then in the definition of the
exponential map we can require ${\mathbf{f}}$ to be \ a smooth immersion. In
this case, the notion of the linearized Poisson structure is still well
defined. But the pull back ${\mathbf{f}}^{\ast}\Psi$ does not isomorphic to
the original Poisson structure $\Psi$ in general.
\end{rem}
\end{dfn}

To prove Theorem 5.1, we will use results obtained in Section~4.

\subsection{The transitive Lie algebroid of a symplectic leaf}

As is well known, the Poisson bracket (5.1) on $M$ admits the natural
extension to the bracket for $1$-forms on $M$:
\begin{equation}
\{\alpha,\beta\}_{T^{\ast}M}=\Psi^{\#}(\alpha)\rfloor d\beta-\Psi^{\#}%
(\beta)\rfloor d\alpha-d\langle\alpha,\Psi^{\#}(\beta)\rangle.
\end{equation}
This structure makes the cotangent bundle $T^{\ast}M$ a Lie algebroid :
\begin{equation}
\left(  T^{\ast}M,\{\,,\,\}_{T^{\ast}M},\rho=\Psi^{\#}\right)
\end{equation}
which is called the Lie algebroid of the Poisson manifold $(M,\Psi)$
\cite{We5}. Notice that if $M$ is not regular, then the Lie algebroid (5.7) is
not transitive.

Given a symplectic leaf $(B,\omega)$ of $M$, one can restrict the bracket
$\{\,,\,\}_{T^{\ast}M}$ to a bracket $\{\,,\,\}_{T_{B}^{\ast}M}$ on smooth
sections of the restricted cotangent bundle $T_{B}^{\ast}M$ . The result is
the transitive Lie algebroid \cite{IKV} ( also see \cite{Ku} for general
criteria of Lie subalgebroids):
\begin{equation}
\left(  T_{B}^{\ast}M,\{\,,\,\}_{T_{B}^{\ast}M},\rho=\rho_{B}\right)
\end{equation}
with anchor
\begin{equation}
\rho_{B}:\,T_{B}^{\ast}M\rightarrow T^{\ast}B\xrightarrow{-(\omega^\flat
)^{-1}}TB,
\end{equation}
where the first morphism is induced by the inclusion $TB\hookrightarrow
T_{B}M$ and $\omega^{\flat}:\,TB\rightarrow T^{\ast}B$ is the bundle map
associated with the symplectic structure $\omega$ $(\omega^{\flat
}(u)=u\,\lrcorner\,\omega)$ The isotropy of this Lie algebroid coincides with
the annihilator $TB^{0}=\ker_{B}\Psi^{\#}$ of $TB$ in $T_{B}M$. We will call
(5.8) the \textbf{transitive Lie algebroid} of the \textbf{symplectic leaf}
$B$.

Let ${\mathcal{N}}$ be the normal bundle to the leaf $B$ and ${\mathbf{f}%
}:\,{\mathcal{N}}\rightarrow M$ be an exponential map. Then the differential
\[
d_{B}{\mathbf{f}}:\,T_{B}{\mathcal{N}}=TB\oplus{\mathcal{N}}\rightarrow T_{B}M
\]
is identical on $TB$ and takes the subbundle ${\mathcal{N}}$ \ to the
complementary subbundle $S=d_{B}{\mathbf{f}}({\mathcal{N}})$ to $TB$. Let
$S^{0}$ be the annihilator of $S$ in $T_{B}M$. The natural splitting
\begin{equation}
T_{B}^{\ast}M=S^{0}\oplus TB^{0}%
\end{equation}
defines the connection $\gamma_{{\mathbf{f}}}:TB\rightarrow T_{B}^{\ast}M$ in
the Lie algebroid (5.8).On the other hand, the exact sequence of vector
bundles $TB\rightarrow T_{B}M\overset{\nu}{\rightarrow}{\mathcal{N}}$
\ induces the dual exact sequence
\begin{equation}
{\mathcal{N}}^{\ast}\overset{\nu^{\ast}}{\rightarrow}T_{B}^{\ast}M\rightarrow
T^{\ast}B\text{.}%
\end{equation}
Using (5.10) and (5.11), we define the vector bundle isomorphism
\begin{equation}
\iota_{{\mathbf{f}}}=\gamma_{{\mathbf{f}}}\oplus\nu^{\ast}:TB\oplus
{\mathcal{N}}^{\ast}\rightarrow T_{B}^{\ast}M
\end{equation}
wich induces the Lie algebroid structure on $A=TB\oplus{\mathcal{N}}^{\ast}:$
\ $\{\,a_{1},\,a_{2}\}_{A}=\iota_{{\mathbf{f}}}^{-1}(\{\iota_{{\mathbf{f}}%
}(a_{1}),\iota_{{\mathbf{f}}}(a_{2})\}_{T_{B}^{\ast}M})$. Thus, we get the
transitive Lie algebroid over $B$ with distinguished connection:
\begin{equation}
\left(  A=TB\oplus{\mathcal{N}}^{\ast},\{\,,\,\}_{A},\rho=\operatorname{pr}%
_{1},\gamma_{0}\right)  \text{.}%
\end{equation}
Here the anchor is the projection onto the first factor, the conormal bundle
${\mathcal{N}}^{\ast}$ is the isotropy and the connection $\gamma_{0}$ is the
canonical injection (4.26) whose pull back via $\iota_{{\mathbf{f}}}$
coincides with the ${\mathbf{f}}$-dependent connection, $\gamma_{{\mathbf{f}}%
}=\iota_{{\mathbf{f}}}\circ\gamma_{0}$.

\bigskip\ Now we can proceed to the proof of Theorem~5.1. Given an exponential
map ${\mathbf{f}}$, we defined the coupling tensor $\Pi^{A,\gamma_{0}}$ on the
normal bundle ${\mathcal{N}}$ associated with the transitive Lie algebroid
$A_{{}}$ in (5.13) and connection $\gamma_{0}$. Clearly $\Pi^{A,\gamma_{0}}$
is equivalent to the coupling tensor associated with the transitive Lie
algebroid of $B$ (5.8) and the connection $\gamma_{{\mathbf{f}}}$. Finally, we
obseve that $\Pi^{A,\gamma_{\mathbf{0}}}$ is just the first approximation to
$\Psi$ at $B$ generated by the exponential map ${\mathbf{f}}$,
\begin{equation}
\Pi^{{\mathbf{f}}}=\Pi^{A,\gamma_{\mathbf{0}}}.
\end{equation}
Here we use the following equivalent reformulation of Definition~5.1.: a
coupling tensor $\Pi$ with an exponential map ${\mathbf{f}}$ defines a first
approximation to $\Psi$ at $B$ if the geometric data of $\Pi$ are obtained
from the geometric data of ${\mathbf{f}}^{\ast}\Psi$ by means of the
linearization at $B$. The independence of $\Pi^{{\mathbf{f}}}$ of the choice
of ${{\mathbf{f}}}$ (up to a nieghborhood equivalence) follows from Theorem~4.2.

We can conclude: \textit{the linearized Poisson structure of }$\Psi
$\textit{\ at a closed symplectic leaf }$(B,\omega)$\textit{\ coincides with
the }$\omega$\textit{-coupling structure of the transitive Lie algebroid of
the leaf.}

Now it is natural to say that a Poisson stucture $\Psi$ is
\textit{linearizable} at a closed symplectic leaf $(B,\omega)$ if there exists
an exponential map ${{\mathbf{f}}}$ such that the pull back ${\mathbf{f}%
}^{\ast}\Psi$ and the first approximation $\Pi^{{\mathbf{f}}}$ are isomorphic
over the zero section $B\hookrightarrow{\mathcal{N}}$ . This definition does
not depend on the choice of ${{\mathbf{f}}}$ .

\begin{rem}
If $\Lambda=0$, then one can try to introduce second approximations to $\Psi$
at $B$, using, for example, results \cite{Du}.
\end{rem}

To end this section, as a consequence of the above results, we give an
affirmative answer to the question on the Poisson realization of transitive
Lie algebroids.

\begin{thm}
Every transitive Lie algebroid $A$ over a connected symplectic base
$(B,\omega)$ can be realized as the transitive Lie algebroid of the symplectic
leaf $(B,\omega)$ of a certain Poisson manifold.
\end{thm}

\subsection{Homotopy invariants}

The notion of the \textit{reduced linear Poisson holonomy} of a symplectic
leaf $B$\textit{,} introduced in \cite{GiGo} (also see \cite{Fe} ), can be
defined as a homotopy invariant of the transitive Lie algebroid of $B$ (5.8).
To see that, pick two coonections $\gamma$ and $\widetilde{\gamma}$ in the Lie
algebroid $T_{B}^{\ast}M$ and consider the corresponding adjoint connections
$\nabla^{\gamma}$ and $\nabla^{\widetilde{\gamma}}$on the isotropy
$\mathcal{L}=TB^{0}$. Fix a point $b_{0}\in B$ and consider a smooth path
$[0$,$1]\ni t\mapsto\sigma(t)\in B$ starting at $b_{0}$, $\sigma(0)=b_{0}$.
Denote by $\mathcal{P}_{t}:\mathcal{L}_{b_{0}}\rightarrow\mathcal{L}%
_{\sigma(t)\text{ }}$and $\widetilde{\mathcal{P}_{t}}:\mathcal{L}_{b_{0}%
}\rightarrow\mathcal{L}_{\sigma(t)\text{ }}$ parallel transport operators
associated with linear connections $\nabla^{\gamma}$ and $\nabla
^{\widetilde{\gamma}}$, respectevely. Define a time dependent field of linear
operators on the fiber $\mathcal{L}_{b_{0}}$ as follows
\begin{equation}
\Xi_{t}:=\mathcal{P}_{t}^{-1}\circ(\operatorname{ad}\circ\mu(\frac{d\sigma
(t)}{dt}))\circ\mathcal{P}_{t}\text{,}%
\end{equation}
where $\mu$ is a $\mathcal{L}$-valued 1-form on $B$, defined in (4.19). It
follows from (4.6) that $\Xi_{t}\in\operatorname{ad}(\mathcal{L}_{b_{0}%
})\approx\operatorname{ad}(\frak{g})$( the adjoint algebra of the typical
fiber $\frak{g}$ ) for all $t\in\lbrack0,1]$. Consider the evolution operator
$\mathbb{T}_{t}\in\operatorname{Ad}(\mathcal{L}_{b_{0}})\approx
\operatorname{Ad}(\frak{g}):$%
\begin{equation}
\frac{d\mathbb{T}_{t}}{dt}=\Xi_{t}\circ\mathbb{T}_{t},\text{ \ \ \ \ \ }%
\mathbb{T}_{0}=\operatorname{id}\text{.}%
\end{equation}

Then we get the following relationship between parallel transports of two
adjoint connections \cite{KV1} : $\widetilde{\mathcal{P}_{t}}=\mathcal{P}%
_{t}\circ\mathbb{T}_{t}$. This implies that for every loop $\sigma\in
\Omega(B;b_{0})$ based at $b_{0}$, the corresponding elements of holonomy
groups $\widetilde{\mathcal{P}_{\sigma}}\in\operatorname{Hol}_{_{b_{0}}%
}^{\nabla^{\widetilde{\gamma}}}\subset\operatorname{Aut}(\frak{g)}$ and
$\mathcal{P}_{\sigma}\in\operatorname{Hol}_{_{b_{0}}}^{\nabla^{\gamma}}%
\subset\operatorname{Aut}(\frak{g)}$ are related by $\widetilde{\mathcal{P}%
_{\sigma}}=\mathcal{P}_{\sigma}\circ\mathbb{T}_{\sigma}$, where $\mathbb{T}%
_{\sigma}\in\operatorname{Ad}(\frak{g})$. Thus, there is a well defined
homorphism $\Omega(B;b_{0})\rightarrow\operatorname{Aut}(\frak{g)/}%
\operatorname{Ad}(\frak{g})$, which does not depend on the choice of an
adjoint connection. If we consider the conjugate homorphism $\Omega
(B;b_{0})\rightarrow\operatorname{Aut}(\frak{g}^{\ast}\frak{)/}%
\operatorname{Ad}(\frak{g}^{\ast})$, then its cotangent lift coincides with
the definition of the reduced linear Poisson holonomy of $B$ given in
\cite{GiGo,Fe}.

\setcounter{section}{0} \renewcommand{\thesection}{\Alph{section}}

\section{Appendix: the proof of Lemma 3.1}

First, remark that if $\Gamma$ is an Ehresmann connection on a fiber bundle
$\pi:\,E\rightarrow B$, then the horizontal lift and the covariant exterior
derivative (2.14) satisfy the modified Cartan formula
\begin{equation}
L_{\operatorname{hor}(u)}=\imath_{u}\circ\partial_{\Gamma}+\partial_{\Gamma
}\circ\imath_{u},\qquad u\in{\mathcal{X}}(B).
\end{equation}
Here $\imath_{u}$ is the interior product. Moreover, the commutator of the
horizontal lift $\operatorname{hor}(u)$ with an arbitrary vertical vector
field is again a vertical vector field,
\begin{equation}
\lbrack\operatorname{hor}(u),{\mathcal{X}}_{V}(E)]\in{\mathcal{X}}_{V}(E).
\end{equation}

Let $\Pi_{t}$ be the time-dependent coupling tensor associated with geometric
data $(\Gamma_{t},{\mathcal{V}},\mathbb{F}_{t})$ in (3.14), (3.15), and let
$X_{t}^{h}\in{\mathcal{X}}_{H}(\mathcal{E}_{0})$ be an arbitrary
time-dependent horizontal vector field. Using properties (A.1), (A.2) and the
standard properties of the Schouten bracket, from relations (2.16)--(2.19) for
$(\Gamma_{t},{\mathcal{V}},\mathbb{F}_{t})$ we deduce the key formula
\begin{align}
L_{X_{t}^{h}}\Pi_{t} &  =-\frac{1}{2}H^{ii^{\prime}}H^{jj^{\prime}}%
(\partial_{\Gamma_{t}}(X_{t}\rfloor\mathbb{F}_{t}))_{i^{\prime}j^{\prime}%
}\operatorname{hor}_{t}(\partial_{i})\wedge\operatorname{hor}_{t}(\partial
_{j})\nonumber\\
&  \qquad+H^{is}({\mathcal{V}}^{\#}d\mathbb{F}_{t}(X_{t},\partial
_{s}))^{\sigma}\partial_{\sigma}\wedge\operatorname{hor}_{t}(\partial_{i}).
\end{align}
Here $\operatorname{hor}_{t}$ is the horizontal lift associated with
$\Gamma_{t}$ and we use the local representations
\[
\Pi_{t}=\frac{1}{2}H^{ij}\operatorname{hor}_{t}(\partial_{i})\wedge
\operatorname{hor}_{t}(\partial_{j})+\frac{1}{2}{\mathcal{V}}^{\sigma
\sigma^{\prime}}\partial_{\sigma}\wedge\partial_{\sigma^{\prime}},
\]
where $\partial_{i}=\partial/\partial\xi^{i}$, $\partial_{\sigma}%
=\partial/\partial x^{\sigma}$, $(\xi^{i})$ and $(x^{\sigma})$ are local
coordinates on the base and the fiber of $\pi$, respectively. Let
$\mathbb{F}_{t}=\frac{1}{2}F_{ij}d\xi^{i}\wedge d\xi^{j}$. Taking into account
$H^{is}F_{sj}=-\delta_{j}^{i^{\prime}}$ and relations (3.14), (3.15), we get
also
\begin{align}
\frac{\partial}{\partial t}\Pi_{t} &  =-\frac{1}{2}H^{ii^{\prime}%
}H^{jj^{\prime}}\frac{\partial}{\partial t}F_{i^{\prime}j^{\prime}%
}\operatorname{hor}_{t}(\partial_{i})\wedge\operatorname{hor}_{t}(\partial
_{j})\nonumber\\
&  \qquad-H^{is}({\mathcal{V}}^{\#}d\phi(\partial_{s}))^{\sigma}%
\partial_{\sigma}\wedge\operatorname{hor}_{t}(\partial_{i}).
\end{align}

Now a direct consequence of (A.3) and (A.4) is that a time-dependent
horizontal vector field $X_{t}^{h}$ is a solution of the homological equation
(3.20) if and only if the associated element $X_{t}\in{\mathcal{X}}(B)\otimes
C^{\infty}(\mathcal{E}_{0})$ satisfies the following two equations
\begin{gather}
\partial_{\Gamma_{t}}(X_{t}\rfloor\mathbb{F}_{t})+\frac{\partial}{\partial
t}\mathbb{F}_{t}=0,\\
X_{t}\rfloor\mathbb{F}_{t}=\phi+c,
\end{gather}
where $c\in\Omega^{1}(B)\otimes\operatorname{Casim}_{{\mathcal{V}}%
}(\mathcal{E}_{0})$ is arbitrary. Taking $c=0$ and $X_{t}$ as the solution of
(3.16), we reduce (A.5) to the identity
\[
\partial_{\Gamma_{t}}\phi=\partial_{\Gamma}\phi+t\{\phi\wedge\phi
\}_{{\mathcal{V}}},\qquad t\in\lbrack0,1],
\]
which holds because of the assumption (3.9). This completes the proof.

\end{document}